\documentclass[12pt,reqno]{amsart}

\usepackage{amscd}
\usepackage{amsfonts}
\usepackage{amsmath}
\usepackage{amssymb}
\usepackage{amsthm}
\usepackage{fancyhdr}
\usepackage{latexsym}
\usepackage[colorlinks=true, pdfstartview=FitV, linkcolor=blue,
            citecolor=blue, urlcolor=blue]{hyperref}

\input epsf
\input texdraw
\input txdtools.tex
\input xy
\xyoption{all}

\usepackage{color}

\definecolor{Red}{rgb}{1.00, 0.00, 0.00}
\definecolor{DarkGreen}{rgb}{0.00, 1.00, 0.00}
\definecolor{Blue}{rgb}{0.00, 0.00, 1.00}
\definecolor{Cyan}{rgb}{0.00, 1.00, 1.00}
\definecolor{Magenta}{rgb}{1.00, 0.00, 1.00}
\definecolor{DeepSkyBlue}{rgb}{0.00, 0.75, 1.00}
\definecolor{DarkGreen}{rgb}{0.00, 0.39, 0.00}
\definecolor{SpringGreen}{rgb}{0.00, 1.00, 0.50}
\definecolor{DarkOrange}{rgb}{1.00, 0.55, 0.00}
\definecolor{OrangeRed}{rgb}{1.00, 0.27, 0.00}
\definecolor{DeepPink}{rgb}{1.00, 0.08, 0.57}
\definecolor{DarkViolet}{rgb}{0.58, 0.00, 0.82}
\definecolor{SaddleBrown}{rgb}{0.54, 0.27, 0.07}
\definecolor{Black}{rgb}{0.00, 0.00, 0.00}
\definecolor{dark-magenta}{rgb}{.5,0,.5}
\definecolor{myblack}{rgb}{0,0,0}
\definecolor{darkgray}{gray}{0.5}
\definecolor{lightgray}{gray}{0.75}

\usepackage[pdftex]{graphicx}
\usepackage{epstopdf}
 
\newcommand{\medno}{\medskip\noindent}

\newcommand{\nin}{\noindent}

\newcommand{\rr}{\mathbb{R}}
\newcommand{\p}{\partial}

\newcommand{\ee}{\varepsilon}

\def\refer #1\par{\noindent\hangindent=\parindent\hangafter=1 #1\par}

\theoremstyle{plain}  
\newtheorem{theorem}{Theorem}
\newtheorem{proposition}{Proposition}
\newtheorem{lemma}{Lemma}

\theoremstyle{definition}

\begin{document}



\title{Non-uniform dependence on initial data\\ for
 the  CH equation on the line
  }
  
  \author{{\it A. Alexandrou Himonas \& Carlos Kenig}}

\date{September  23, 2010 
[Corrected version of:
A. Himonas and C. Kenig,
{\it Non-uniform dependence on initial data for the CH equation on the line,}
Differential and Integral Equations,
Vol. 22, No. 3-4 (2009) pp. 201-224.]}



\begin{abstract}
For $s>3/2$ two sequences of  CH solutions living in a bounded
set of the Sobolev space $H^s(\rr)$ are constructed,
whose distance at the initial time  is  converging to zero
while at any later time is bounded  below by 
a positive constant.
This implies that the  solution map  of the CH equation
is not uniformly continuous in $H^s(\rr)$.
\end{abstract}

\maketitle
\markboth{ Non-uniform dependence for CH  equation on the line}
                    {Alex Himonas and Carlos Kenig}

\parindent0in
\parskip0.1in


%
%

\section{Introduction}
 \setcounter{equation}{0}
We consider the Cauchy problem for the Camassa-Holm equation (CH)
\begin{equation} 
\label{CH}
\partial_t u + u\partial_x u
+
\partial_x
\Big(
1 - \partial_x^2
\Big)^{-1}
 \Big[
u^2 + \frac{1}{2}(\partial_x u)^2
\Big]
= 0, 
\end{equation}
\begin{equation} 
\label{CH-data} 
u(x, 0) = u_0 (x),\,\,
\ x \in  \mathbb{R},  \;\;  \ t \in \mathbb{R}.
\end{equation}
This equation  appeared initially in the context of
hereditary symmetries studied by Fuchssteiner and Fokas \cite{ff}.
However, it was written explicitly as a water wave equation
by Camassa and Holm \cite{ch},  who showed  that CH is biHamiltonian
and  studied its ``peakon" solutions. Since then CH has been  rederived 
in various ways  by  Misio\l ek \cite{mi},  Johnson \cite{j}, 
Constantin and Lannes  \cite{cl},  and Ionescu-Kruse \cite{i}.

Well-posedness on the line was first established by Li and Olver.
In  \cite{lo} they showed that if $s>3/2$ then CH is locally well-posed in $H^s(\rr)$ 
with solutions depending {\it continuously} on  initial data.
The proof was based on a regularization technique similar to
the one used by Bona and Smith for the KdV equation \cite{bs}.
A similar  result has also been proved by   Rodriguez-Blanco \cite{rb} 
by using Kato's theory  for quasilinear  equations \cite{k}.
Moreover,  global well-posedness in $H^1(\rr)$ for the CH equation
has been studied by Bressan and Constantin in \cite{bc}.
However, well-posednes of CH in $H^s(\rr)$ for $s\in (1, 3/2]$ 
remains  an open question.

In this paper,  we show  that dependence of CH solutions on initial data 
in Sobolev spaces can not be better than continuous.
More precisely,  we  prove  the following result.

\begin{theorem}
\label{CH-non-unif-dependence} 
 If  $s>3/2$  then the  flow  map $u_0 \to u(t)$   
for the CH equation is not uniformly  continuous from
any  bounded set of $H^s(\rr)$ into $C([-T, T]; H^s(\rr))$.
More precisely,   there  exist two sequences of  CH solutions
$u_n(t)$ and $v_n(t)$ in $C([-T, T]; H^s(\rr))$ such that
\begin{equation} 
\label{H-s-bdd} 
\| u_n(t)  \|_{H^s(\rr)}
  +
  \| v_n(t)  \|_{H^s(\rr)}
\lesssim
 1,
\end{equation}
  \begin{equation} 
  \label{zero-limit-at-0}
  \lim_{n\to\infty}
\|
u_n(0)
-
v_n(0)
  \|_{H^s(\rr)}
=
0,
\end{equation}
and
  \begin{equation} 
  \label{bdd-away-from-0}
  \liminf_{n\to\infty}
\|
u_n(t)
-
v_n(t)
  \|_{H^s(\rr)}
\gtrsim
    \sin t, 
    \quad
    |t|<T\le 1.
\end{equation}
 \end{theorem}
For $s=1$  Theorem  \ref{CH-non-unif-dependence} has been
already proved by Himonas, Misio\l ek and Ponce in  \cite{hmp} by using traveling wave solutions  that are smooth except at finitely many points at
which the slope is  $\pm \infty$ (cuspons).
Also, in   \cite{hmp} the analogues result for the periodic CH was proved.
 For  $s\ge 2$ non-uniform  continuity of the CH solution map 
in the periodic case was established in \cite{hm} using
high frequency traveling wave solutions
and  following an approach  similar
to the one used in \cite{kpv} by Kenig, Ponce and Vega.
We mention that this method does not
work in the non-periodic case because 
the traveling wave solutions  do not 
live in $H^s(\rr)$.

Also, it is worth mentioning the following implication
of Theorem  \ref{CH-non-unif-dependence} concerning ways
for proving local well-posedness for CH. 
 The fact that the data-to-solution map is not uniformly continuous
from any  bounded set of $H^s(\rr)$ into $C([-T, T]; H^s(\rr))$
 tells us that  local well-posedness of CH  in $H^s$ cannot be established by
a solely contraction principle argument.

 The proof of Theorem \ref{CH-non-unif-dependence}  is based on the method 
of approximate solutions  used by Koch and Tzvetkov in \cite{kt}
and Christ, Colliander and Tao in \cite{cct}.
The idea is to choose approximate solutions  consisting of
a low-frequency part and a high-frequency part,
which satisfy  the three conclusions of Theorem \ref{CH-non-unif-dependence}.
Furthermore, solving the Cauchy problem with 
initial data given by evaluating the approximate solutions at $t = 0$ 
must yield actual solutions whose difference from 
the approximate solutions is  negligible. 

The literature about CH is extensive.
For some other results  about this equation
we refer the reader  to McKean \cite{mc},
Constantin and Strauss \cite{cs},
Himonas, Misio\l ek, Ponce and Zhou \cite{hmpz},
and  Molinet's  survey article  \cite{mo}.

The paper is structured as follows. 
In section 2 we recall the well-posedness result of Li and Olver and
use it to prove  the basic energy estimate (see \eqref{CH-diff-ineq})
from which  we derive a lower bound for the lifespan of the solution
as well   an estimate of the $H^s$ norm of the solution  $u(t)$
in terms of the  the $H^s$  norm of the initial data $u_0$
(see Proposition \ref{Lifespan-u-size}).
In section 3 we construct  approximate solutions consisting 
of a low-frequance part and a high-frequency part,
 and compute the error. 
In section 4 we estimate the $H^1$-norm of this error.
In  section 5 we solve the Cauchy problem for the CH equation with 
initial data given by the approximate solutions evaluated at time zero,
and estimate the $H^1$-norm of the difference beween actual and 
approximate solutions (see Lemma \ref{CH-differ-H1-est-lem}).
Finally, in section 6 we conclude with the proof of Theorem \ref{CH-non-unif-dependence}.

%
%
%
%

\section{Local well-posedness} 
\setcounter{equation}{0}
          
We  shall need the following well-posedness result,
 proved in \cite{lo}  using 
a regularization technique.
\begin{theorem}
\label{CH-wp}
[Li-Olver]
Suppose that the function $u_0(x)$ belongs to the Sobolev space 
$H^s(\rr)$ for some  $s >3/2$. Then there is  a $T>0$, which depends 
only on  $\|u_0\|_{H^s}$, such that there exists a unique function $u(x, t)$
solving  the Cauchy problem  \eqref{CH}--\eqref{CH-data} 
in the sense of distributions with  $u \in C([0, T]; H^s)$.
When $s\ge 3$,  $u$ is also a classical solution to \eqref{CH}--\eqref{CH-data}. 
Moreover, the solution $u$ depends continuously  on the initial data $u_0$
in the sense that the mapping of the initial data to the solution 
is continuous from the Sobolev space $H^s$ to the space $C([0, T]; H^s)$.
\end{theorem}

Using the information provided by Theorem \ref{CH-wp},
next we shall prove an explicit estimate for the time of existence  $T$
of the solution  $u(t)$. Also, we  will  show that at any time $t$ in the time interval
$[0, T]$ the $H^s$ norm of the solution  $u(t)$ is  dominated
by the  $H^s$ norm of the initial data $u_0$.

\begin{proposition}
\label{Lifespan-u-size}
Let  $s>3/2$.  If  $u$ is  the solution  of  the Cauchy problem  \eqref{CH}--\eqref{CH-data} 
described in Theorem \ref{CH-wp}  then  its lifespan
(the maximal existence time)
 is greater than 
     \begin{equation}
   \label{Lifespan-size}
   T
   \doteq
   \frac{1}{2c_s}
   \frac{1}{ \|u_0 \|_{H^s(\rr)}},
 \end{equation}
where $c_s$  is a constant depending only on $s$.
 Also, we have that
  \begin{equation}
   \label{u-u0-Hs-bound}
\|u(t)\|_{H^s(\rr))}
  \le
  2
  \|u_0 \|_{H^s(\rr)},
  \quad
  0\le t \le T.
   \end{equation}
  \end{proposition}

\noindent
{\bf Proof.}
  The derivation of the lower bound for the   lifespan  \eqref{Lifespan-size} 
and  the solution size estimate \eqref{u-u0-Hs-bound} is based on the following differential
inequality for the solution $u$
\begin{equation} 
\label{CH-diff-ineq}
\frac 12
\frac{d}{dt}
  \|u(t)\|_{H^{s}(\rr)}^2
\le
c_s
 \|u(t)\|_{H^{s}(\rr)}^3,
 \quad
 0\le t \le T.
\end{equation}
This inequality  can be extracted  from the proof  of Theorem \ref{CH-wp}
in \cite{lo} using the  energy estimate (3.6) proved   for the following regularization
$$
\p_tu - \p_x^2\p_tu +\ee\p_x^4\p_tu+3u\p_xu -\p_xu\p_x^2 u-u\p_x^3u=0
$$ 
of the CH equation
and  letting $\ee$ go to zero. Here,  we shall prove  inequality \eqref{CH-diff-ineq} by
following the approach used for quasilinear symmetric
hyperbolic systems  in Taylor  \cite{t1}.

  For any $s\in \rr$ let   $D^s=(1-\p_x^2)^{s/2}$ be the  operator  defined by 
$$
\widehat{D^s f}(\xi) \doteq (1 + \xi^2)^{s/2} \widehat{f}(\xi),
$$
where $ \widehat{f}$ is the Fourier transform
$$
\widehat{f}(\xi) =  \int_\rr e^{-ix\xi}\,\,  f(x) dx.
$$
Then for   $f \in H^s(\rr)$ we have
$$
\|f\|_{H^s(\rr)}
=
  \int_\rr  \big(1 + \xi^2\big)^{s}| \widehat{f}(\xi)|^2 \frac{d\xi}{2\pi}
  =
 \|D^sf\|_{L^2(\rr)}.
$$
Now let $u$ be the solution  to the Cauchy problem   \eqref{CH}--\eqref{CH-data},
which  according to Theorem \ref{CH-wp}  belongs  in $ C([0, T]; H^s)$.
Solving   \eqref{CH} for $\p_t u$ we obtain
\begin{equation} 
\label{CH-ut-form}
\partial_t   u 
=
-u\partial_x u
-
D^{-2}
\partial_x
 \Big[
u^2 + \frac{1}{2}(\partial_x u)^2
\Big].
\end{equation} 
Starting with \eqref{CH-ut-form}  we  want to derive the energy estimate in $H^s$
expressed by inequality \eqref{CH-diff-ineq}.
We can  form  $\frac{d}{dt}\|u\|_{H^s(\rr)}^2$ by applying formally  the operator  $D^s$
 to both sides  of  \eqref{CH-ut-form}, then multiply the resulting equation by  $D^su$
and integrate it with respect to $x$. Note that since $u\in H^s$ the second term in the right-hand side of  \eqref{CH-ut-form}  is in  $H^s$ too. 
However, the first term, that is the product  $u\partial_x u$
is   only in $H^{s-1}$. To deal with this problem we will replace  \eqref{CH-ut-form} 
by its molified  smooth version 
\begin{equation} 
\label{CH-ut-form-moli}
\partial_t J_\ee u 
=
-J_\ee(u\partial_x u) 
-
D^{-2}
\partial_x
 \Big[
J_\ee(u^2) + \frac{1}{2}J_\ee[(\partial_x u)^2]
\Big],
\end{equation} 

where  for each $\ee\in (0, 1]$ the operator $J_\ee$ is  the Friedrichs mollifier 
defined by  
\begin{equation} 
\label{Fried-moli}
J_\ee f(x)= j_\ee \ast f (x).
\end{equation}
Here  $j(x)$ is  a  $C^\infty$ function supported in the interval $[-1, 1]$ such that
$j(x)\ge 0$,  $\int_\rr j(x)dx =1$
and 
$$
j_\ee(x)= \frac{1}{\ee}  \, j\big(\frac{x}{\ee}\big).
$$
Applying  the operator $D^s$ to  both sides of  \eqref{CH-ut-form-moli},
then  multiplying the resulting equation by $D^s J_\ee u$
and integrating it for $x\in\rr$ gives
\begin{equation} 
\label{CH-moli-int}
\begin{split}
\frac 12
\frac{d}{dt} \|J_\ee u \|_{H^s}^2
&=
-
   \int_\rr
     D^sJ_\ee(u\partial_x u) \cdot   D^sJ_\ee u   \, dx
- 
       \int_\rr
   D^{s-2} \partial_x  J_\ee(u^2)  \cdot   D^sJ_\ee u \, dx
     \\
&-
 \frac 12
      \int_\rr
      D^{s-2} \partial_x J_\ee[(\partial_x u)^2]  \cdot   D^sJ_\ee u  \, dx.
      \end{split}
\end{equation}
In what follows next we use the fact that  $D^s$ and $J_\ee$ commute
and  that  $J_\ee$ satisfies the properties 
\begin{equation} 
\label{J-e-inner-prod-property}
(J_\ee f, g)_{L^2}=( f, J_\ee g)_{L^2},
\end{equation}
and
\begin{equation} 
\label{Je-u-Hs}
    \| J_\ee u \|_{H^s}
  \le
       \|  u \|_{H^s}.
\end{equation}
%

  %
  %
  %
  \noindent
{\bf Estimating the Burgers term.}
To estimate the  first integral in the right-hand side of \eqref{CH-moli-int}
we write it as follows
\begin{equation} 
\label{int1-est-calc1}
\begin{split}
   \int_\rr
     D^sJ_\ee(u\partial_x u) \cdot   D^sJ_\ee u   \, dx
 &= 
  \int_\rr
     D^s(u\partial_x u) \cdot   J_\ee D^sJ_\ee u   \, dx
 \\
  &=
    \int_\rr
\big[ 
D^s(u\p_x u)  -  u D^s (\p_xu)
\big]
  J_\ee D^sJ_\ee u   \, dx
 \\
  &+
       \int_\rr
  u D^s (\p_xu)
\cdot   J_\ee D^sJ_\ee u   \, dx.
\end{split}
\end{equation}
Now, we estimate the first term  in the  right-hand side of  \eqref{int1-est-calc1}.
Applying the Cauchy-Schwarz inequality gives
\begin{equation} 
\label{int1-est-calc2}
\begin{split}
\Big|
    \int_\rr
\big[ 
D^s(u\p_x u)  -  u D^s (\p_xu)
\big]
 J_\ee D^sJ_\ee u   \, dx
\Big|
&\le
\|
D^s(u\p_x u)  -  u D^s (\p_xu)
\|_{L^2}
\|
J_\ee D^sJ_\ee u
\|_{L^2}
\\
&\le
\|
D^s(u\p_x u)  -  u D^s (\p_xu)
\|_{L^2}
\|
 u
\|_{H^s}
\\
&\le
 2c_s    \| \p_x u \|_{L^\infty} 
    \| u \|_{H^s}^2,
\end{split}
\end{equation}
where the last step follows from the estimate  
\begin{equation} 
\label{int1-est-calc3}
 \|  D^s(u\p_x u)  -  u D^s (\p_xu) \|_{L^2}
 \le
 2c_s    \| \p_x u \|_{L^\infty} 
    \| u \|_{H^s},
\end{equation}
which we prove below by using   the  following Kato-Ponce commutator 
estimate \cite{kp} (see also  Ionescu and  Kenig  \cite{ik}).
 
 \begin{lemma} 
 \label{KP-lemma}
 [Kato-Ponce]
  If  $s>0$ then there is $c_s>0$ such that  for any $f, g \in H^s(\rr)$ 
 \begin{equation} 
\label{KP-com-est}
 \| D^{s} \big(fg) -  f D^s g\|_{L^2}
 \le
 c_s\big(
    \| D^{s}f \|_{L^2}    \| g \|_{L^\infty} 
    +
     \| \p_xf \|_{L^\infty}    \| D^{s-1}g \|_{L^2}   
 \big).
 \end{equation}
 \end{lemma}
 In fact, applying  this estimate with $f=u$ and $g=\p_xu$ gives 
 \begin{equation} 
\label{int1-est-calc4}
\begin{split}
  \|  D^s(u\p_x u)  -  u D^s (\p_xu) \|_{L^2}
&
 \le
 c_s\big(
    \| D^{s}u \|_{L^2}    \| \p_x u \|_{L^\infty} 
    +
     \| \p_xu \|_{L^\infty}    \| D^{s-1}\p_x u \|_{L^2}   
 \big)
 \\
 & \le
 c_s    \| \p_x u \|_{L^\infty} 
 \big(
    \| D^{s}u \|_{L^2}  
    +
  \| D^{s}u \|_{L^2}   
 \big) 
  \\
 & \le
 2c_s    \| \p_x u \|_{L^\infty} 
    \| u \|_{H^s}, 
\end{split}
  \end{equation}
 which  is the desired estimate  \eqref{int1-est-calc3}.

 Next, we  estimate the  second  integral in the right-hand side of \eqref{int1-est-calc1}.
 Note  if there were no $J_\ee$'s  involved  then this would have been  done in
  a straightforward manner as follows
 \begin{equation} 
 \label{int1-est-calc5}
 \begin{split}
 \Big|
        \int_\rr
  u D^s (\p_xu)
\cdot  D^su   \, dx
 \Big|
 &=
  \Big|
       \frac 12 \  \int_\rr
  u \p_x\big[(D^s u)^2\big]  \, dx
   \Big|
   \\
  & =
    \Big|
      - \frac 12 \  \int_\rr
   \p_xu \, (D^s u)^2 \, dx
     \Big|
        \\
        &\le
\frac 12 \| \p_x u \|_{L^\infty} 
    \| u \|_{H^s}^2.
   \end{split}
  \end{equation}
When the $J_\ee$'s are  involved  the idea is the same.  However, the implementation is 
 more technical since we need  to commute $J_\ee$ so that  is
 grouped correctly. We accomplish this as follows 
\begin{equation} 
\label{int1-est-calc6}
\begin{split}
    \int_\rr
  u D^s (\p_xu)
\cdot   J_\ee D^sJ_\ee u   \, dx
 &= 
    \int_\rr
  J_\ee  u D^s (\p_xu)
\cdot  D^sJ_\ee u   \, dx
  \\
  &=
    \int_\rr
    \Big(
  [J_\ee,  u ] D^s (\p_xu) +    u J_\ee  D^s (\p_xu)
  \Big)
\cdot  D^sJ_\ee u   \, dx
 \\
  &=
    \int_\rr
  [J_\ee,  u ]\p_x D^s u 
  \cdot  D^sJ_\ee u   \, dx
   \\
  &+
   \int_\rr
u  \p_xD^s J_\ee u
\cdot  D^sJ_\ee u   \, dx.
\end{split}
\end{equation}
Estimating the  second integral  of  the right-hand  side of \eqref{int1-est-calc6}
 like  we have done in \eqref{int1-est-calc5} we get 
 \begin{equation} 
 \label{int1-est-calc7}
 \begin{split}
 \Big|
          \int_\rr
u  \p_xD^s J_\ee u
\cdot  D^sJ_\ee u   \, dx
 \Big|
 &=
  \Big|
       \frac 12 \  \int_\rr
  u \p_x\big[(D^s J_\ee u)^2\big]  \, dx
   \Big|
   \\
  & =
    \Big|
      - \frac 12 \  \int_\rr
   \p_xu \, (D^s J_\ee u)^2 \, dx
     \Big|
        \\
        &\le
\frac 12 \| \p_x u \|_{L^\infty} 
    \| J_\ee u \|_{H^s}^2
       \\
        &\le
\frac 12 \| \p_x u \|_{L^\infty} 
    \|  u \|_{H^s}^2.
   \end{split}
  \end{equation}
For estimating   the first  integral  of  the right-hand  side of \eqref{int1-est-calc6}
 we apply the Cauchy-Schwarz inequality and  we have
\begin{equation} 
\label{int1-est-calc8}
  \begin{split}
    \Big|
    \int_\rr
  [J_\ee,  u ]\p_x D^s u 
  \cdot  D^sJ_\ee u   \, dx
      \Big|
   &\le
    \|    [J_\ee,  u ]\p_x D^s u    \|_{L^2}
       \|   D^sJ_\ee u   \|_{L^2}
          \\
        &\le
          \|    [J_\ee,  u ]\p_x D^s u    \|_{L^2}
       \|  u   \|_{H^s}
               \\
        &\le
        c
      \| \p_x u \|_{L^\infty} 
    \|  u \|_{H^s}^2,
 \end{split}
      \end{equation}
where the last step of the above inequality is justified by the following result.
\begin{lemma}
\label{Je-u-com}  Let $u(x)$ be a function such that    $\| \p_x u \|_{L^\infty} <\infty$.
Then, there is $c>0$ such that for   any $f \in L^2(\rr)$ we have
\begin{equation} 
\label{Je-u-com-L2-est}
 \|    [J_\ee,  u ]\p_x f  \|_{L^2}
     \le
        c
      \| \p_x u \|_{L^\infty} 
    \|  f \|_{L^2}.
      \end{equation}
\end{lemma}
{\bf Proof.} We have 
\begin{equation} 
\label{Je-u-com-calc1}
\begin{split}
  [J_\ee,  u ]\p_x f  (x)
     &=
  J_\ee( u \p_x f )(x)   -   u J_\ee(\p_x f )(x)
  \\
  &=
  j_\ee  \ast ( u \p_x f )(x) 
  -
  u(x)   (j_\ee  \ast  \p_x f) (x) 
    \\
  &=
  \int_\rr
    j_\ee(x-y) u(y) f'(y) \, dy
    -
    u(x)  \int_\rr
    j_\ee(x-y)f'(y) \, dy
     \\
  &=
  \int_\rr
  \frac{1}{\ee}
    j\big(\frac{x-y}{\ee}\big) \big[u(y)-u(x)\big] f'(y) \, dy.
      \end{split}
      \end{equation}
Integrating by parts  and using the mean value theorem gives
\begin{equation} 
\label{Je-u-com-calc2}
\begin{split}
  [J_\ee,  u ]\p_x f  (x)
     &=
  -
 \int_{\rr}
  \frac{1}{\ee}
    j\big(\frac{x-y}{\ee}\big) u'(y) f(y) \, dy
       \\
    &
  +
 \int_{\rr}
  \frac{1}{\ee^2}
    j' \big(\frac{x-y}{\ee}\big)  \big[u(y)-u(x)\big] f(y) \, dy
    \\
    &=
  -
 \int_{|y-x|<\ee}
  \frac{1}{\ee}
    j\big(\frac{x-y}{\ee}\big) u'(y) f(y) \, dy
    \\
    &+
 \int_{|y-x|<\ee}
  \frac{1}{\ee^2}
    j' \big(\frac{x-y}{\ee}\big) u'(\xi(x, y))(y-x) f(y) \, dy.
 \end{split}
      \end{equation}
Above we have used our assumption  that  $j(x)$ is supported on the interval $[-1, 1]$. So, 
using the bound  $|(x-y)/\ee| <1$ and taking  absolute values 
we obtain that
\begin{equation} 
\label{Je-u-com-calc3}
\begin{split}
  \big| [J_\ee,  u ]\p_x f  (x)  \big| 
     &\le
         \| \p_x u \|_{L^\infty} 
         \Big(
 \int_{\rr}
  \frac{1}{\ee}
    j\big(\frac{x-y}{\ee}\big) |f(y)| \, dy
       \\
    &
 +
 \int_{\rr}
  \frac{1}{\ee}
   \big|   j' \big(\frac{x-y}{\ee}\big) \big|  |f(y)| \, dy  
    \Big)
     \\
    &=
           \| \p_x u \|_{L^\infty} 
         \Big(
j_\ee \ast |f| (x)
 +
|j_\ee' |\ast |f| (x)
    \Big).
     \end{split}
      \end{equation}
Finally, applying Young's inequality we get 
\begin{equation} 
\label{Je-u-com-calc4}
\begin{split}
 \|    [J_\ee,  u ]\p_x f  \|_{L^2}
  &\le
     \| \p_x u \|_{L^\infty} 
      \big(
     \|  j_\ee  \|_{L^1}      \|  f \|_{L^2}
     +
         \|  j_\ee'  \|_{L^1}      \|  f \|_{L^2}
         \big)
         \\
        & =
      \big(
     \|  j  \|_{L^1} +   \|  j' \|_{L^1} 
      \big)  
            \| \p_x u \|_{L^\infty}  
         \|  f \|_{L^2},
  \end{split}
      \end{equation}
      which gives  the desired inequality \eqref{Je-u-com-L2-est}
      with constant $c=  \|  j  \|_{L^1} +   \|  j' \|_{L^1}$. 
      \hfil $\square$
  Combining the  inequalities  \eqref{int1-est-calc1}, \eqref{int1-est-calc2},
\eqref{int1-est-calc7} and  \eqref{int1-est-calc7}  we obtain the
following estimate for the Burgers term of the CH equation
\begin{equation} 
\label{Burgers-energy-est}
\Big|
   \int_\rr
     D^sJ_\ee(u\partial_x u) \cdot   D^sJ_\ee u   \, dx
     \Big|
     \le
     c_s
      \| \p_x u \|_{L^\infty} 
    \|  u \|_{H^s}^2.
\end{equation}
%

  %
  %
  %
  \noindent
{\bf Estimating the nonlocal   $D^{s-2} \partial_x  J_\ee(u^2)$.}
To  estimate the  second  integral in the right-hand side of \eqref{CH-moli-int}
  we apply  the Cauchy-Schwarz inequality  and we get
\begin{equation} 
\label{int2-est-calc1}
\begin{split}
\Big|
    \int_\rr
   D^{s-2} \partial_x  J_\ee(u^2)  \cdot   D^sJ_\ee u \, dx
   \Big|
 &\le
     \|     D^{s-2} \partial_x  J_\ee(u^2)    \|_{L^2}
         \|      D^sJ_\ee u  \|_{L^2}
         \\
         &\le
              \|   u^2 \|_{H^{s-1}}     \|  u  \|_{H^s}
            \\
         &\le
              \|   u^2 \|_{H^s}    \|  u  \|_{H^s}.
           \end{split}
\end{equation}
Now, we use the following  estimate for  the Sobolev norm of a product,
which can be found in Taylor \cite{t2}  (see Corollary 10.6).
For any $s>0$ and $1<p<\infty$  there is $C=C_{s,p}>0$ such that
\begin{equation} 
\label{best-Sob-product-est}
\|  f g \|_{H^{s,p}}                
         \le
   C  \Big[
                 \|  f \|_{H^{s,p}}    \|  g\|_{L^\infty}     +   \|  f \|_{L^\infty}     \|  g\|_{H^{s,p}}
                 \Big]. 
           \end{equation}
Using this result with $s=2$ and $f=g=u$  from 
\eqref{int2-est-calc1}  we obtain that
\begin{equation} 
\label{nonloc-u2-energy-Sob-best}
\Big|
    \int_\rr
   D^{s-2} \partial_x  J_\ee(u^2)  \cdot   D^sJ_\ee u \, dx
   \Big|
\le
2 c_s 
\|  u \|_{L^\infty}    \|  u\|_{H^s}^2.
\end{equation}
%

  %
  %
  %
    \noindent
{\bf Estimating the nonlocal  term  $ D^{s-2} \partial_x J_\ee[(\partial_x u)^2] $.}
 As before,
applying  the Cauchy-Schwarz inequality  we have
\begin{equation} 
\label{nonloc-ux2-energy-Sob-best}
\begin{split}
\Big|
  \int_\rr
      D^{s-2} \partial_x J_\ee[(\partial_x u)^2]  \cdot   D^sJ_\ee u  \, dx
   \Big|
 &\le
     \|  D^{s-2} \partial_x J_\ee[(\partial_x u)^2]  \|_{L^2}
         \|      D^sJ_\ee u  \|_{L^2}
         \\
         &\le
              \|  (\p_x u)^2 \|_{H^{s-1}}     \|  u  \|_{H^s}
            \\
         &\le
         c_s
              \|   \p_x u \|_{H^{s-1}}^2    \|  u  \|_{H^s}
               \\
         &  \le
              2 c_s       \|  \p_xu  \|_{L^\infty}     \|  u  \|_{H^s}^2,
            \end{split}
\end{equation}
where in the last step we used estimate \eqref{best-Sob-product-est}
applied with   $s$ replace by $s-1>0$ and $f=g=\p_xu$.

Now, combining equation
 \eqref{Fried-moli}  and  estimates 
 \eqref{Burgers-energy-est}, 
  \eqref{nonloc-u2-energy-Sob-best}, 
   \eqref{nonloc-ux2-energy-Sob-best}
  we obtain the differential inequality
\begin{equation} 
\label{CH-moly-ineq}
\frac 12
\frac{d}{dt}
  \|J_\ee u(t)\|_{H^s}^2
\le
c_s
 \| u(t) \|_{C^1} 
 \|u(t)\|_{H^s}^2,
 \quad
 0\le t \le T.
\end{equation}
Next, integrating \eqref{CH-moly-ineq}  from 0 to $t$, $t<T$, gives
\begin{equation} 
\label{CH-moly-ineq-int1}
\frac 12
  \|J_\ee u(t)\|_{H^s}^2
  -
  \frac 12
  \|J_\ee u(0)\|_{H^s}^2 
\le
c_s
\int_0^t
 \| u(\tau) \|_{C^1} 
 \|u(\tau)\|_{H^s}^2\, d\tau.
\end{equation}
Then, letting $\ee$ go to $0$  \eqref{CH-moly-ineq-int1}  gives 
\begin{equation} 
\label{CH-moly-ineq-int2}
 \frac 12
  \|u(t)\|_{H^s}^2
  -
  \frac 12
  \|u(0)\|_{H^s}^2  
\le
c_s
\int_0^t
 \| u(\tau) \|_{C^1} 
 \|u(\tau)\|_{H^s}^2\, d\tau.
\end{equation}
Finally,  from \eqref{CH-moly-ineq-int2}  using Gronwall's inequality  we obtain
the following lemma, which summarizes  our estimates thus far.

\begin{lemma}
\label{CH-energy-inequality}
 Let   $s>3/2$ and   $u \in C([0, T]; H^s)$  be   the solution
of the  Cauchy problem  \eqref{CH}--\eqref{CH-data}.  Then
\begin{equation} 
\label{CH-energy-ineq}
\frac 12
\frac{d}{dt}
  \|u(t)\|_{H^s}^2
\le
c_s
 \| u(t) \|_{C^1} 
 \|u(t)\|_{H^s}^2,
 \quad
 0\le t \le T.
\end{equation}
\end{lemma}
Since $s>3/2$ using Sobolev's inequality 
\begin{equation} 
\label{Sob-C1-ineq}
 \| u(t) \|_{C^1}
 \le
c_s
 \|u(t)\|_{H^s}, 
\end{equation}
from \eqref{CH-energy-ineq} we obtain the desired inequality \eqref{CH-diff-ineq}.

%

%
%
%
%
%
%
\vskip0.1in
\noindent
{\bf  Lifespan estimate.}    
To derive an explicit formula for  $T=T( \|v(0)\|_{H^{s}})$ we proceed as follows.
Letting  $y(t)=  \|u(t)\|_{H^{s}}^2$   inequality  \eqref{CH-diff-ineq}  
takes the form
\begin{equation} 
\label{energy-y-ineq}
\frac 12
 y^{-3/2}\frac{dy}{dt}
\le 
c_s,
  \qquad
  y(0)=y_0=  \|u_0\|_{H^{s}}^2.
\end{equation}
Integrating  \eqref{energy-y-ineq} from  0  to $t$ gives
\begin{equation} 
\label{energy-y-ineq-calc1}
\frac{1}{\sqrt{y_0}}  - \frac{1}{\sqrt{y(t)}} 
\le 
c_s t.
\end{equation}
Replacing $y(t)$ with   $\|u(t)\|_{H^{s}}^2$  and solving for  $\|u(t)\|_{H^s}$
we obtain the formula
\begin{equation} 
\label{norm-u(t)-formula}
 \|u(t)\|_{H^s}
\le
\frac{ \|u_0\|_{H^s}}{1-c_s\|u_0\|_{H^s} t}.
\end{equation}
Now,  from \eqref{norm-u(t)-formula} we see that  $\|u(t)\|_{H^s}^2$ is finite  if 
\begin{equation*} 
\label{Lifespan-calc1}
c_s    \|u_0\|_{H^s} t<1,
\end{equation*}
or
\begin{equation} 
\label{Lifespan-calc1}
t
<
\frac{1}{ c_s \|u_0\|_{H^s}}.
\end{equation}
Therefore, the  solution  $u(t)$ to the CH   Cauchy problem certainly  exists  
for $0\le t <T_0$, where
\begin{equation} 
\label{CH-Lifespan}
T_0
=
\frac{1}{ c_s \|u_0\|_{H^s}}.
\end{equation}
%

%
%
%
%
%
%
\noindent
{\bf  Size of the solution estimate.} If we choose  $T=1/2 T_0$, that is
\begin{equation} 
\label{T-def}
T
=
\frac{1}{2 c_s \|u_0\|_{H^s}},
\end{equation}
then  for $0\le t \le T$ inequality  \eqref{norm-u(t)-formula} gives 
\begin{equation*} 
\label{u(t)-u(0)-bound}
\|u(t)\|_{H^{s}}
\le
\frac{ \|u_0\|_{H^s}}{1-(c_s\|u_0\|_{H^s})/(2 c_s \|u_0\|_{H^s})},
\end{equation*}
or 
\begin{equation} 
\label{u(t)-u(0)-bound}
\|u(t)\|_{H^{s}}
\le
2  \|u_0\|_{H^s},
\quad 
0\le t \le T.
\end{equation}
This completes the proof  of Proposition  \ref{Lifespan-u-size}.
\hfil $\square$

%
%
%
%
%
%
\vskip0.1in
\noindent
{\bf   Remark.}   Inequality   \eqref{CH-energy-ineq} can be used to show that 
 if  $u \in C([0, T]; H^s)$, $s>3/2$,  is   a solution
of  Cauchy problem  \eqref{CH}--\eqref{CH-data} 
 such that
$
\sup_{0\le t<T}
 \| u(t) \|_{C^1}
 <\infty
$
then  $u(t)$ persists to be a solution beyond the time $T$.
In particular,  we can show that if  the lifespan  $T$ of $u$
is finite then  
$
\sup_{0\le t<T}
 \| u(t) \|_{C^1} = \infty
$
(see   Theorem 6.2 in \cite{lo}).

%
%
%
%
%
%
%
%
%
%
%
%
%
\section{Construction  of approximate solutions }
 \setcounter{equation}{0}
 Here we shall construct  a two-parameter family of approximate solutions
 $u^{\omega, \lambda}=u^{\omega, \lambda}(x, t)$,
 each member of which  consists of two parts, that is
\begin{equation} 
\label{approx-sln} 
u^{\omega, \lambda}
=
u_\ell + u^h. 
\end{equation}
The high frequency  part  $u^h$  is  given by
\begin{equation} 
\label{high-frequency-approx-sln}
u^h=u^{h,\omega, \lambda}(x, t)
 = 
  \lambda^{-\delta/2  -  s}
  \varphi (\frac{x}{\lambda^{\delta}})
  \cos (\lambda x-\omega t),
 \end{equation}
and  is not a solution of CH.  
Here   $\varphi$ is a $C^{\infty}$ function  such that
\begin{equation} 
\label{phi}
\varphi(x)=
\begin{cases}
&1,  \text { if  } |x|<1,\\
&0,  \text { if  } |x|\ge 2.
\end{cases}
 \end{equation}
 The low frequency  part  
 $
 u_\ell =u_{\ell, \omega, \lambda}(x, t)
 $
  is the  solution to the  following Cauchy problem for  CH
 \begin{equation} 
\label{low-frequency-sln}
\partial_t u_\ell + u_\ell\partial_x u_\ell
+
\Lambda^{-1}
 \Big[
u_\ell^2 + \frac{1}{2}(\partial_x u_\ell)^2
\Big]
= 0, 
\end{equation}
\begin{equation} 
\label{low-frequency-data} 
u_\ell(x, 0) 
=
 \omega  \lambda^{-1}
\tilde{\varphi}
(\frac{x}{\lambda^{\delta}}), \ x \in  \mathbb{R},  \;\;  \ t \in \mathbb{R},
\end{equation}
where  $\tilde{\varphi}$  is a $C_0^\infty(\mathbb{R})$ function 
 such that 
\begin{equation} 
\label{phi-relation}
\tilde{\varphi}(x)=1,
 \,\, \text{if} \,\,
 x\in \text{supp }\varphi.
 \end{equation}
Furthermore,
 $\Lambda^{-1}$ denotes the order $-1$ pseudodifferential operator
 \begin{equation} 
\label{Lambda-1-def}
\Lambda^{-1}
=
\partial_x
\Big(
1 - \partial_x^2
\Big)^{-1}.
\end{equation}
As it is explained in Lemma  \ref{u-low-cp-info} below,
the initial value problem  
\eqref{low-frequency-sln}--\eqref{low-frequency-data}
has a unique 
smooth solution $u_\ell$ belonging
in $H^s(\rr)$ for all $s$. Thus, the approximate solutions
$u^{\omega, \lambda}$  belong in every Sobolev space.

Substituting the approximate solution 
$
u^{\omega, \lambda}
=
u_\ell + u^h 
$
into  CH  equation  we obtain the following expression
\begin{equation*} 
\label{approx-sln-error} 
\begin{split}
F
&=
\p_t
 u^h
+
u_\ell 
\p_x u^h
+ 
u^h
\p_x u_\ell
+ 
u^h
\p_x u^h
+
\Lambda^{-1}
\Big[
2u_\ell u^h
+
(u^h)^2
+
\p_xu_\ell \p_xu^h 
+
\frac 12 (\p_x u^h)^2
\Big]
\\
&+
\p_t u_\ell + u_\ell\p_x u_\ell
+
\partial_x
\Big(
1 - \partial_x^2
\Big)^{-1}
 \Big[
u_\ell^2 + \frac{1}{2}(\partial_x u_\ell)^2
\Big].
\end{split}
\end{equation*}
Now,  taking into consideration that $u_\ell$ solves 
CH we obtain the following error for the approximate solution 
\begin{equation} 
\label{approx-sln-error} 
\begin{split}
F
&=
\p_t
 u^h
+
u_\ell 
\p_x u^h
+ 
u^h
\p_x u_\ell
+ 
u^h
\p_x u^h
\\
&
+
\Lambda^{-1}
\Big[
2u_\ell u^h
+
(u^h)^2
+
\p_xu_\ell \p_xu^h 
+
\frac 12 (\p_x u^h)^2
\Big].
\end{split}
\end{equation}
Computing $\p_t u^h$  gives 
\begin{equation} 
\label{t-der-of-u-high-1}
\p_tu^h(x, t)
=
 \omega
  \lambda^{-\delta/2  -  s}
  \varphi (\frac{x}{\lambda^{\delta}})
  \sin (\lambda x-\omega t ).
 \end{equation}
 
 Furthermore, since $\tilde{\varphi}$ is equal to 1 on the support of $\varphi$
we see that  we can write $\p_t u^h$ in the following form
\begin{equation} 
\label{t-der-of-u-high}
\begin{split}
\p_tu^h(x, t)
 &= 
 \omega
   \tilde{\varphi} (\frac{x}{\lambda^{\delta}})
  \lambda^{-\delta/2  -  s}
  \varphi (\frac{x}{\lambda^{\delta}})
  \sin (\lambda x-\omega t )
\\
&=
\lambda 
u_\ell(x, 0)
\cdot 
 \lambda^{-\delta/2  -  s}
  \varphi (\frac{x}{\lambda^{\delta}})
  \sin (\lambda x-\omega t ).
\end{split}
 \end{equation}
 Computing the spacial derivative of $u^h$ gives
\begin{equation} 
\label{x-der-of-u-high}
\begin{split}
\p_xu^h(x, t)
 &= 
-
 \lambda 
\cdot 
  \lambda^{-\delta/2  -  s}
  \varphi (\frac{x}{\lambda^{\delta}})
  \sin (\lambda x-\omega t),
\\
&+
  \lambda^{-\frac 32 \delta  -  s}
  \p_x\varphi (\frac{x}{\lambda^{\delta}})
  \cos (\lambda x-\omega t ).
\end{split}
 \end{equation}
Then, using  \eqref{t-der-of-u-high} and \eqref{x-der-of-u-high}  
we find that  
\begin{equation} 
\label{dt-uh-calc} 
\begin{split}
\p_t
 u^h
+
u_\ell 
\p_x u^h
&=
\lambda
\Big[
u_\ell(x, 0) - u_\ell(x, t) 
\Big]
 \lambda^{-\delta/2  -  s}
  \varphi (\frac{x}{\lambda^{\delta}})
  \sin (\lambda x-\omega t )
  \\
  &+
  u_\ell(x, t) 
\cdot
    \lambda^{-\frac 32 \delta  -  s}
  \p_x\varphi (\frac{x}{\lambda^{\delta}})
  \cos (\lambda x-\omega t ).
\end{split}
\end{equation}
Therefore, the error \eqref{approx-sln-error}  of the approximate solution $u^{\omega, \lambda}$ 
is given by
\begin{equation} 
\label{F-sum-of Fj} 
F
=
F_1 + F_2+ \cdots +F_8,
\end{equation}
where
\begin{equation} 
\label{F-j-def} 
\begin{split}
&F_1
=
\lambda
\Big[
u_\ell(x, 0) - u_\ell(x, t) 
\Big]
 \lambda^{-\delta/2  -  s}
  \varphi (\frac{x}{\lambda^{\delta}})
  \sin (\lambda x-\omega t )
  \\
  &
  F_2=
  u_\ell(x, t) 
\cdot
    \lambda^{-\frac 32 \delta  -  s}
  \p_x\varphi (\frac{x}{\lambda^{\delta}})
  \cos (\lambda x-\omega t )
   \\
  &
  F_3   =   u^h \p_x u_\ell
  \\
  &
  F_4=  u^h  \p_x u^h
    \\
  &
  F_5=  \Lambda^{-1} \big[  2u_\ell u^h \big]
   \\
  &
  F_6=     \Lambda^{-1} \big[    (u^h)^2  \big]
    \\
  &
  F_7=    \Lambda^{-1} \big[    \p_xu_\ell \p_xu^h    \big]
    \\
  &
  F_8=    \Lambda^{-1} \big[   \frac 12 (\p_x u^h)^2     \big].
\end{split}
\end{equation}

Next we shall estimate the size of the error  $F$.

%
%
%
%

\section{  Estimating the $H^1$ norm of the error } 
\setcounter{equation}{0}

To estimate the $H^1$ norm  of the error  $F$
it suffices to estimate  the $H^1$ norm  of  each term $F_j$. 
Observe that  each $F_j$ is expressed in terms of 
$u_\ell$ and $u^h$. The high frequency part   $u^h$
is defined by  formula \eqref{high-frequency-approx-sln}
and   
     \begin{equation}
   \label{Hs-norm-of-u-h}
  \|
  u^h(t)
  \|_{H^s(\rr)}
\approx
1,
\quad 
\text{for }
\,
\lambda>>1,
\end{equation}
because of the following result.

%
%
%
%
%
%
\begin{lemma}
   \label{lem:Hs-norm-of-ap-sl}
   Let 
  $\psi \in \mathcal{S}(\rr)$, $1<\delta<2$ and   $\alpha  \in \rr$.
Then for any  $s\ge 0$  we have that 
     \begin{equation}
   \label{Hs-norm-of-ap-sl-2}
   \lim_{\lambda\to \infty}
   \lambda^{-\frac{1}{2}\delta-s}
  \|
  \psi (\frac{x}{\lambda^\delta}) \cos (\lambda x - \alpha)
  \|_{H^s(\rr)}
=
\frac{1}{\sqrt{2}}
 \|
  \psi
  \|_{L^2(\rr)}.
 \end{equation}
Relation \eqref{Hs-norm-of-ap-sl-2} is also true if  $\cos$ is replaced by   $\sin$.
\end{lemma}
Although this lemma   can be found in \cite{kt}, we include its proof here for
the convenience  of the reader.

\vskip0.1in
\noindent
{\bf Proof.}  Since
 \begin{equation*}
 \begin{split}
\Big( 
\psi (\frac{x}{\lambda^\delta})  \cos (\lambda x - \alpha)
\widehat {\Big)}(\xi)
=
\frac 12
\lambda^\delta
\big[
e^{-i\alpha}
\widehat{
\psi}
(\lambda^\delta (\xi -\lambda))
+
e^{i\alpha}
\widehat{
\psi}
(\lambda^\delta (\xi + \lambda))
\big],
 \end{split}
 \end{equation*}
 we have that 
 \begin{equation*}
 \label{Lambda-1-trig-calc}
 \begin{split}
    \lambda^{-\delta-2s}
    \|
\psi(\frac{x}{\lambda^{\delta}})
 \cos (\lambda x-\alpha)
  \|_{H^s(\rr)}^2
&
=
\frac{ \lambda^{-2s+\delta}}{8\pi}
\int_\rr
\big(1+\xi^2)^s
\big|
e^{-i\alpha}
\widehat{
\psi}
(\lambda^\delta (\xi -\lambda))
+
e^{i\alpha}
\widehat{
\psi}
(\lambda^\delta (\xi + \lambda))
\big|^2
d\xi
\\
&=
   \frac{ \lambda^{-2s+\delta}}{8\pi}
   \Big[
\int_\rr
\big(1+\xi^2)^s
\big|
\widehat{
\psi}
(\lambda^\delta (\xi -\lambda))
\big|^2
d\xi
\\
&+
\int_\rr
\big(1+\xi^2)^s
\big|
\widehat{
\psi}
(\lambda^\delta (\xi +\lambda))
\big|^2
d\xi
\\
&+
2
\int_\rr
\big(1+\xi^2)^s
\text{Re}
 \big[
e^{-2i\alpha}
\widehat{
\psi}
(\lambda^\delta (\xi - \lambda))
\bar{
\widehat{\psi}
}
(\lambda^\delta (\xi + \lambda))
\big]
d\xi
\Big].
 \end{split}
 \end{equation*}
Now,  in the first and third integral we make 
the change of variables $\eta=\lambda^\delta (\xi-\lambda)$,
while in the second we let  $\eta=\lambda^\delta (\xi+\lambda)$.
Thus, we have
 \begin{equation*}
 \label{Lambda-1-trig-calc}
 \begin{split}
    \lambda^{-\delta-2s}
    \|
\psi(\frac{x}{\lambda^{\delta}})
 \cos (\lambda x-\alpha)
  \|_{H^s(\rr)}^2
&
=
\frac{ \lambda^{-2s}}{8\pi}
   \Big[
\int_\rr
\Big(1+\big(\frac{\eta}{\lambda^\delta}+\lambda)^2
\Big)^s
\big|
\widehat{
\psi}
(\eta)
\big|^2
d\eta
\\
&+
\int_\rr
\Big(1+\big(\frac{\eta}{\lambda^\delta}- \lambda)^2
\Big)^s
\big|
\widehat{
\psi}
(\eta)
\big|^2
d\eta
\\
&+
2
\int_\rr
\Big(1+\big(\frac{\eta}{\lambda^\delta}+\lambda)^2
\Big)^s
\text{Re}
 \big[
e^{-2i\alpha}
\widehat{
\psi}
(\eta)
\bar{
\widehat{\psi}
}
(\eta+2\lambda^{\delta+1})
\big]
d\xi
\Big].
 \end{split}
 \end{equation*}
 Moving the factor $\lambda^{-2s}$ inside the integrals gives
 \begin{equation*}
 \label{Lambda-1-trig-calc}
 \begin{split}
    \lambda^{-\delta-2s}
    \|
\psi(\frac{x}{\lambda^{\delta}})
 \cos (\lambda x-\alpha)
  \|_{H^s(\rr)}^2
&
=
\frac{1}{8\pi}
   \Big[
\int_\rr
\Big(
\frac{1}{\lambda^2}+\big(\frac{\eta}{\lambda^{\delta+1}}+ 1)^2
\Big)^s
\big|
\widehat{
\psi}
(\eta)
\big|^2
d\eta
\\
&+
\int_\rr
\Big(
\frac{1}{\lambda^2}+\big(\frac{\eta}{\lambda^{\delta+1}} - 1)^2
\Big)^s
\big|
\widehat{
\psi}
(\eta)
\big|^2
d\eta
\\
&+
2
\int_\rr
\Big(
\frac{1}{\lambda^2}+\big(\frac{\eta}{\lambda^{\delta+1}}+ 1)^2
\Big)^s
\text{Re}
 \big[
e^{-2i\alpha}
\widehat{
\psi}
(\eta)
\bar{
\widehat{\psi}
}
(\eta+2\lambda^{\delta+1})
\big]
d\xi
\Big].
 \end{split}
 \end{equation*}
 Since   $\psi \in \mathcal{S}(\rr)$ we have that
  $\widehat{\psi} (\eta+2\lambda^{\delta+1}) \to 0$ as $\lambda \to \infty$.
  Therefore, applying the  dominated convergence theorem  we see that
 the third  integral goes to zero  while each of the other two goes  to $\|\widehat{\psi}\|_{L^2}^2$.
 Therefore, we obtain that 
 \begin{equation*}
 \label{Lambda-1-trig-calc}
 \lim_{\lambda \to \infty}
    \lambda^{-\delta-2s}
    \|
\psi(\frac{x}{\lambda^{\delta}})
 \cos (\lambda x-\alpha)
  \|_{H^s(\rr)}^2
=
\frac{1}{4\pi}
\|\widehat{\psi}\|_{L^2}^2
=
\frac{1}{2}
\|\psi\|_{L^2}^2,
 \end{equation*}
 which proves the lemma. \,\, $\square$
 
 %
 %
 %
%
%
%
%

 %
 %

As we have stated earlier, the low frequency part  $u_\ell$
is  the solution of the Cauchy problem
\eqref{low-frequency-sln}--\eqref{low-frequency-data}.
Next lemma summarizes the basic information  
about  $u_\ell$.

\begin{lemma}
\label{u-low-cp-info}
Let  $\omega$ be bounded,
$0<\delta <2$
 and $\lambda>>1$.
Then,   the  initial value problem 
 \eqref{low-frequency-sln}--\eqref{low-frequency-data}
 has  a unique   smooth solution  $u_\ell
 \in  C([0, 1]; H^s(\rr))$,  for all $s > 3/2$,  and satisfying the 
 estimate
     \begin{equation}
        \label{Hs-norm-of u-ell-t-est}
     \|
    u_\ell(t)
  \|_{H^s(\rr)}
\le
c_s
 \lambda^{-1+\delta/2},
 \quad
 0\le t \le 1.
     \end{equation}
     \end{lemma}
{\bf Proof.}  Let  $s\ge 0$.
   For any    function $\psi\in \mathcal{S}(\rr)$ we have
     \begin{equation}
        \label{psi-delta-estmate}
     \|
\psi(\frac{x}{\lambda^{\delta}})
  \|_{H^s(\rr)}
\le
        \lambda^{\delta/2}
        \,\,
         \|
\psi
       \|_{H^s(\rr)}.
     \end{equation}
    In fact,  using the relation  
 $\widehat{\psi(x/\rho)} (\xi)
 =
 \rho
 \widehat{\psi}(\rho \xi) 
 $
 and making the change of variables $\eta=\lambda^\delta \xi$ 
 we obtain
 \begin{equation*}
   \label{psi-delta-est}
   \begin{split}
  \|
\psi(\frac{x}{\lambda^{\delta}})
\|_{H^s}^2
&=
\frac{1}{2\pi}
   \int_{\rr}
  (1+\xi^2)^s
   |
\lambda^\delta
\,
\widehat{
\psi}
(\lambda^\delta \xi)
|^2
  d\xi 
\\
  &
  =
\frac{1}{2\pi}
   \int_{\rr}
\Big(1+\frac{\eta^2}{\lambda^{2\delta}}\Big)^s
   \cdot
\lambda^{2\delta}
|
\widehat{
\psi}
(\eta)
|^2
\,
  \frac{d\eta}{\lambda^\delta}
  \\
  &
  =
\lambda^{\delta} 
\cdot
 \frac{1}{2\pi}
   \int_{\rr}
\Big(1+\frac{\eta^2}{\lambda^{2\delta}}\Big)^s
|
\widehat{
\psi}
(\eta)
|^2
\,
  d\eta
  \\
  &
  \le
\lambda^{\delta} 
\cdot
 \frac{1}{2\pi}
   \int_{\rr}
   \big(1+ \eta^2  \big)^s
|
\widehat{
\psi}
(\eta)
|^2
\,
 d\eta
   \\
  &
     =
        \lambda^{\delta}
        \,\,
         \|
\psi
  \|^2_{H^s(\rr)}.
 \end{split}
 \end{equation*}
Now, using  inequality  \eqref{psi-delta-estmate} we  have that
the initial data  $u_\ell (0)$ satisfy the estimate
     \begin{equation}
        \label{Hs-norm-of u-ell-0-est}
     \|
    u_\ell(0)
  \|_{H^s(\rr)}
\le
|\omega|
 \lambda^{-1+\delta/2}
        \,\,
         \|
\tilde{\varphi}
  \|_{H^s(\rr)},
     \end{equation}
     which for $\omega$ bounded decays  if 
     \begin{equation}
        \label{delta-condition}
  \delta < 2.
     \end{equation}

Next,  using  estimate \eqref{Lifespan-size} from
Proposition \ref{Lifespan-u-size}   we have
that the lifespan  $T$  of the solution $u_\ell(t)$ satisfies 
     \begin{equation*}
        \label{Lifespan-est}
        T
        \ge
        \frac{1}{ 2c_s
     \|
    u_\ell(0)
  \|_{H^s(\rr)}}
   \ge
        \frac{c_s'}{ 
   \lambda^{-1+\delta/2}} \ge 1,
   \quad 
   \text{ for }
   \lambda  >>1,
    \end{equation*}
     since $\delta<2$.
  Finally,  if $s \ge 0$ then  from  estimate \eqref{u-u0-Hs-bound}
     of Proposition \ref{Lifespan-u-size}   we  have
     \begin{equation*}
        \label{u-ell-t-est}
     \|
    u_\ell(t)
  \|_{H^s(\rr)}
  \le
       \|
    u_\ell(t)
  \|_{H^{s+2}(\rr)}
  \le
   c_s
        \|
    u_\ell(0)
  \|_{H^{s+2}(\rr)}
\le
   c_s
   \lambda^{-1+\delta/2}.
   \quad \square
     \end{equation*}

       \vskip0.1in    
Now we are ready to estimate  the $H^1$ norm of each error $F_j$.

%
%
%
\medno
   {\bf Estimating the  $H^1$-norm of $F_1$.}   We have 
 \begin{equation}
 \label{H1-est-F1-calc1}
 \begin{split}
   \|F_1\|_{H^1(\rr)}
&=
  \|
 \lambda
\Big[
u_\ell(x, 0) - u_\ell(x, t) 
\Big]
 \lambda^{-\delta/2  -  s}
  \varphi (\frac{x}{\lambda^{\delta}})
  \sin (\lambda x-\omega t )\|_{H^1(\rr)}
\\
&=
 \lambda^{1-\delta/2  -  s}
  \|
    \varphi (\frac{x}{\lambda^{\delta}})
      \sin (\lambda x-\omega t )
  \big[
u_\ell(x, 0) - u_\ell(x, t) 
\big]
\|_{H^1(\rr)}.
\end{split}
 \end{equation}
 Using  the  inequality 
  \begin{equation}
 \label{H1-product}
   \|
f g
\|_{H^1(\rr)}
\le 
\sqrt{2}
\,
   \|
 f
\|_{C^1(\rr)}
   \|
 g
\|_{H^1(\rr)},
 \end{equation}
from \eqref{H1-est-F1-calc1}  we get
\begin{equation*}
 \label{H1-est-F1-calc2}
 \begin{split}
   \|F_1\|_{H^1(\rr)}
&\lesssim
 \lambda^{1-\delta/2  -  s}
  \|
    \varphi (\frac{x}{\lambda^{\delta}})
      \sin (\lambda x-\omega t )
      \|_{C^1(\rr)}
      \|
u_\ell(x, 0) - u_\ell(x, t) 
\|_{H^1(\rr)}
\end{split}
 \end{equation*}
 And, since 
 $\|
    \varphi (\frac{x}{\lambda^{\delta}})
      \sin (\lambda x-\omega t )
      \|_{C^1(\rr)}
      =
      \|  \varphi\|_{L^\infty} \lambda
      $
      the last inequality gives
 \begin{equation}
 \label{H1-est-F1-calc2}
   \|F_1\|_{H^1(\rr)}
\lesssim
 \lambda^{2-\delta/2  -  s}
  \|
u_\ell(x, 0) - u_\ell(x, t) 
\|_{H^1(\rr)}.
 \end{equation}
  To estimate the $H^1$ norm of the difference  $u_\ell(t)-u_\ell(0)$ we apply
 the fundamental theorem of calculus in the time variable to obtain
 \begin{equation}
 \label{FTC-in-t}
 u_\ell(x, t) - u_\ell(x, 0)=\int_0^t \p_t u_\ell (x, \tau) d\tau.
 \end{equation}
 Then, taking the $H^1$ norm  of the space variable to both sides
 of  \eqref{FTC-in-t}
 and passing the norm inside the integral gives
 \begin{equation}
 \label{FTC-in-t-H1}
   \|
u_\ell(x, 0) - u_\ell(x, t) 
\|_{H^1(\rr)}
\le
\int_0^t
   \| \p_t u_\ell (x, \tau)
   \|_{H^1(\rr)}
    d\tau,
     \qquad
 t\in [0, 1].
 \end{equation}
 Next we estimate 
 $ 
  \| \p_t u_\ell (x, \tau)
   \|_{H^1(\rr)}.
$
For this  we solve  equation  \eqref{low-frequency-sln} for  $\partial_t u_\ell$
 to get
 \begin{equation} 
\label{t-der-u-low-relation}
\p_t u_\ell (x, \tau)
=
-
 u_\ell\partial_x u_\ell
-
\Lambda^{-1}
 \big[
u_\ell^2 + \frac{1}{2}(\partial_x u_\ell)^2
\big].
\end{equation}
Thus, at any time in $[0, T]$ we have 
 \begin{equation} 
\label{t-der-u-low-H1}
 \| \p_t u_\ell (x, \tau)
   \|_{H^1(\rr)}
\le
\|
 u_\ell\partial_x u_\ell
 \|_{H^1(\rr)}
+
\|
\Lambda^{-1}
 \big[
u_\ell^2 + \frac{1}{2}(\p_x u_\ell)^2
\big]
 \|_{H^1(\rr)}
\end{equation}
Now, using the inequality
 \begin{equation} 
\label{H1-algebra}
  \| 
fg
  \|_{H^1 (\rr)}
\le c
   \|  
   f 
  \|_{H^1(\rr)}
    \|  
g
  \|_{H^1(\rr)}
\end{equation}
and the estimate 
 \begin{equation} 
\label{Lambda-1-est}
  \| 
  \Lambda^{-1} 
f
  \|_{H^1 (\rr)}
\le
   \|  
f 
  \|_{L^2(\rr)}
\end{equation}
from  \eqref{t-der-u-low-H1} we obtain  that
 \begin{equation} 
 \begin{split}
\label{t-der-u-low-H1-2}
 \| \p_t u_\ell (x, \tau)
   \|_{H^1(\rr)}
   &\lesssim
    \|  
 u_\ell    
  \|_{H^1 (\rr)}
      \|  
 \p_xu_\ell    
  \|_{H^1 (\rr)}
+
\|
u_\ell^2 + \frac{1}{2}(\partial_x u_\ell)^2
  \|_{L^2 (\rr)}
    \\
  &\lesssim
 \|  
 u_\ell    
  \|_{H^1 (\rr)}
      \|  
 u_\ell    
  \|_{H^2 (\rr)}
+
\|
u_\ell^2 
  \|_{L^2 (\rr)}
  +
  \|
   (\partial_x u_\ell)^2
  \|_{L^2 (\rr)}    
  \\
  &\lesssim
      \|  
 u_\ell    
  \|_{H^2 (\rr)}^2
+
\|
u_\ell
  \|_{L^\infty (\rr)}
  \|
u_\ell 
  \|_{L^2 (\rr)}
  +
    \|
   \partial_x u_\ell
  \|_{L^\infty(\rr)}
  \|
   \partial_x u_\ell
  \|_{L^2 (\rr)}
    \\
  &\lesssim
      \|  
 u_\ell    
  \|_{H^2 (\rr)}^2
+
      \|  
 u_\ell    
  \|_{H^1 (\rr)}^2
  +
         \|  
 u_\ell    
  \|_{H^2 (\rr)}^2
      \\
  &\lesssim
      \|  
 u_\ell    
  \|_{H^2 (\rr)}^2.
        \end{split}
\end{equation}
   Using estimate \eqref{Hs-norm-of u-ell-t-est}, from the last inequality we get
 \begin{equation} 
\label{t-der-u-low-H1-fin}
 \| \p_t u_\ell (x, \tau)
   \|_{H^1(\rr)}
\lesssim
  \lambda^{-2+\delta}.
  \end{equation}
  Substituting \eqref{t-der-u-low-H1-fin} into  \eqref{FTC-in-t-H1} we obtain
 \begin{equation}
 \label{H1-of-u-ell-dif-est}
   \|
u_\ell(x, 0) - u_\ell(x, t) 
\|_{H^1(\rr)}
\lesssim
  \lambda^{-2+\delta}.
 \end{equation}
  Finally, combining  \eqref{H1-of-u-ell-dif-est} and  \eqref{H1-est-F1-calc2}
  gives 
 \begin{equation}
 \label{H1-est-F1-calc4}
   \|F_1\|_{H^1(\rr)}
\lesssim
 \lambda^{2-\delta/2  -  s}
 \cdot
 \lambda^{-2+\delta},
 \end{equation}
 which gives 
 \begin{equation}
 \label{H1-est-F1}
   \|F_1\|_{H^1(\rr)}
\lesssim
 \lambda^{- s + \delta/2},
 \qquad
 \lambda>>1.
 \end{equation}
 %
 %

%
%
%
%
\medno
   {\bf Estimating the  $H^1$-norm of $F_2$.}   Reading $F_2$ from  \eqref {F-j-def}  we have
 \begin{equation}
 \label{H1-est-F2-calc1}
 \begin{split}
   \|F_2\|_{H^1(\rr)}
&=
  \|
  u_\ell(x, t) 
\cdot
    \lambda^{-\frac 32 \delta  -  s}
  \p_x\varphi (\frac{x}{\lambda^{\delta}})
  \cos (\lambda x-\omega t )
  \|_{H^1(\rr)}
\\
&\lesssim
 \lambda^{-\frac 32 \delta  -  s}
  \|
   \p_x\varphi (\frac{x}{\lambda^{\delta}})
  \cos (\lambda x-\omega t )
  \|_{C^1(\rr)}
  \|
  u_\ell(x, t) 
\|_{H^1(\rr)}
\\
&\lesssim
 \lambda^{-\frac 32 \delta  -  s}
  \|
   \p_x\varphi (\frac{x}{\lambda^{\delta}})
  \cos (\lambda x-\omega t )
  \|_{C^1(\rr)}
  \|
  u_\ell(x, t) 
\|_{H^2(\rr)}
\\
&\lesssim
 \lambda^{-\frac 32 \delta  -  s}
\cdot 
\lambda
\cdot
 \lambda^{-1+ \frac 12 \delta},
 \end{split}
 \end{equation}
 which gives
 \begin{equation}
 \label{H1-est-F2}
   \|F_2\|_{H^1(\rr)}
\lesssim
 \lambda^{ - s - \delta}.
 \end{equation}
 %
 %

%
%
%
%
\medno
   {\bf Estimating the  $H^1$-norm of $F_3$.}   From  \eqref {F-j-def}  we have
 \begin{equation*}
 \label{H1-est-F3-calc1}
 \begin{split}
   \|F_3(t)\|_{H^1(\rr)}
&=
  \|
  u^h(t)
  \p_xu_\ell(t)
  \|_{H^1(\rr)}
\\
&\lesssim
  \|
  u^h(t)
  \|_{C^1(\rr)}
  \|
 \p_xu_\ell(t)
\|_{H^1(\rr)}
\\
&\lesssim
   \|
  u^h(t)
  \|_{C^1(\rr)}
  \|
 u_\ell(t)
\|_{H^2(\rr)}.
 \end{split}
 \end{equation*}
Using  formula  \eqref{high-frequency-approx-sln} for $u^h$ 
and estimate    \eqref{Hs-norm-of u-ell-t-est} for  $u_\ell$,
from  the last inequality we obtain that
 \begin{equation*}
 \label{H1-est-F3-calc2}
   \|F_3(t)\|_{H^1(\rr)}
\lesssim
 \lambda^{-\frac 12 \delta  -  s +1}
\cdot
 \lambda^{-1+ \frac 12 \delta},
 \end{equation*}
which gives
 \begin{equation}
 \label{H1-est-F3}
   \|F_3(t)\|_{H^1(\rr)}
\lesssim
 \lambda^{ - s},
 \qquad 
 \lambda>>1.
 \end{equation}
 %
 %

%
%
%
%
\medno
   {\bf Estimating the  $H^1$-norm of $F_4$.}   Reading $F_4$ from  \eqref {F-j-def}  
   and using \eqref{best-Sob-product-est} we have
 \begin{equation}
 \label{H1-est-F4-calc1}
 \begin{split}
   \|F_4(t)\|_{H^1(\rr)}
&=
  \|
  u^h(t)
  \p_xu^h(t)
  \|_{H^1(\rr)}
\\
&\lesssim
  \|
  u^h(t)
  \|_{H^1(\rr)}
  \|
 \p_xu^h(t)
\|_{L^\infty(\rr)}
+
  \|
  u^h(t)
  \|_{L^\infty(\rr)}
  \|
 \p_xu^h(t)
  \|_{H^1(\rr)}
\\
&\lesssim
  \|
  u^h(t)
  \|_{H^1(\rr)}
  \|
 \p_xu^h(t)
\|_{L^\infty(\rr)}
+
  \|
  u^h(t)
  \|_{L^\infty(\rr)}
  \|
 u^h(t)
  \|_{H^2(\rr)}.
 \end{split}
 \end{equation}
 Since 
 $$
   \|
  u^h(t)
  \|_{L^\infty(\rr)}
  \lesssim
  \lambda^{ -\frac 12 \delta - s},
  \quad
     \|
  \p_xu^h(t)
  \|_{L^\infty(\rr)}
  \lesssim
  \lambda^{ -\frac 12 \delta - s+1},
 $$
 and  since, by  Lemma  \ref{lem:Hs-norm-of-ap-sl}, we have
\begin{equation*}
 \label{Hk-u-h}
 \begin{split}
     \|
 u^h(t)
\|_{H^k(\rr)}
&=
 \lambda^{-\delta/2  -  s}
     \|
  \varphi (\frac{x}{\lambda^{\delta}})
  \cos (\lambda x-\omega t )
  \|_{H^k(\rr)}
\\
&=
 \lambda^{-s  +  k}
 \cdot
 \lambda^{-\delta/2  - k}
     \|
  \varphi (\frac{x}{\lambda^{\delta}})
  \cos (\lambda x-\omega t )
  \|_{H^k(\rr)}
\\
&\lesssim
   \lambda^{-s  +  k},
 \end{split}
 \end{equation*}
  estimate \eqref{H1-est-F4-calc1} gives
 \begin{equation*}
 \label{H1-est-F4-calc2}
   \|F_4(t)\|_{H^1(\rr)}
\lesssim
  \lambda^{-s  +  1}
  \cdot
   \lambda^{ -\frac 12 \delta - s+1}
+
   \lambda^{ -\frac 12 \delta - s}
\cdot
  \lambda^{-s  +  2}.
 \end{equation*}
Thus,
 \begin{equation}
 \label{H1-est-F4}
   \|F_4(t)\|_{H^1(\rr)}
\lesssim
 \lambda^{ -2 s -\frac 12 \delta +2},
 \qquad 
 \lambda>>1.
 \end{equation}
 %
 %

%
%
%
%
\medno
   {\bf Estimating the  $H^1$-norm of $F_5$.}  We have
 \begin{equation*}
 \label{H1-est-F5-calc1}
 \begin{split}
   \|F_5\|_{H^1(\rr)}
&=
  \|
  \Lambda^{-1} \big[  2u_\ell u^h \big]
  \|_{H^1(\rr)}
\\
&\le
2
  \|
u_\ell u^h
\|_{L^2(\rr)}
\\
&\lesssim
   \|
  u^h
  \|_{L^\infty(\rr)}
  \|
 u_\ell
\|_{L^2(\rr)}
\\
&\lesssim
   \|
  u^h
  \|_{L^\infty(\rr)}
  \|
 u_\ell
\|_{H^2(\rr)}
\\
&\lesssim
   \lambda^{ -\frac 12 \delta - s}
\cdot
  \lambda^{-1  +  \frac 12 \delta },
 \end{split}
 \end{equation*}
 which gives
 \begin{equation}
 \label{H1-est-F5}
   \|F_5(t)\|_{H^1(\rr)}
\lesssim
 \lambda^{ -s -1},
 \qquad 
 \lambda>>1.
 \end{equation}
 %
 %

%
%
%
%
\medno
   {\bf Estimating the  $H^1$-norm of $F_6$.}  From  \eqref {F-j-def}  
   and Lemma  \ref{lem:Hs-norm-of-ap-sl} we have
 \begin{equation*}
 \label{H1-est-F6-calc1}
 \begin{split}
   \|F_5\|_{H^1(\rr)}
&=
  \|
  \Lambda^{-1} \big[    (u^h)^2  \big]
  \|_{H^1(\rr)}
\\
&\le
  \|
(u^h)^2 
\|_{L^2(\rr)}
\\
&\lesssim
   \|
  u^h
  \|_{L^\infty(\rr)}
  \|
  u^h
\|_{L^2(\rr)}
\\
&\lesssim
   \lambda^{ -\frac 12 \delta - s}
\cdot
  \lambda^{-s },
 \end{split}
 \end{equation*}
 which gives
 \begin{equation}
 \label{H1-est-F6}
   \|F_6(t)\|_{H^1(\rr)}
\lesssim
 \lambda^{ -2s - \frac 12 \delta},
 \qquad 
 \lambda>>1.
 \end{equation}
 %
 %

%
%
%
%
\medno
   {\bf Estimating the  $H^1$-norm of $F_7$.}  Also, we  have
 \begin{equation*}
 \label{H1-est-F7-calc1}
 \begin{split}
   \|F_7\|_{H^1(\rr)}
&=
  \|
\Lambda^{-1} \big[    \p_xu_\ell \p_xu^h    \big]
  \|_{H^1(\rr)}
\\
&\le
  \|
\p_xu_\ell \p_xu^h 
\|_{L^2(\rr)}
\\
&\lesssim
   \|
\p_xu^h 
  \|_{L^\infty(\rr)}
  \|
 \p_xu_\ell
\|_{L^2(\rr)}
\\
&\lesssim
   \|
\p_xu^h 
  \|_{L^\infty(\rr)}
  \|
 u_\ell
\|_{H^2(\rr)}
\\
&\lesssim
   \lambda^{ -\frac 12 \delta - s +1}
\cdot
  \lambda^{-1 +\frac 12 \delta },
 \end{split}
 \end{equation*}
 which gives
 \begin{equation}
 \label{H1-est-F6}
   \|F_6(t)\|_{H^1(\rr)}
\lesssim
 \lambda^{ -s },
 \qquad 
 \lambda>>1.
 \end{equation}
 %
 %

%
%
%
%
\medno
   {\bf Estimating the  $H^1$-norm of $F_8$.}   Finally, we  have
 \begin{equation*}
 \label{H1-est-F8-calc1}
 \begin{split}
   \|F_8\|_{H^1(\rr)}
&=
  \|
\Lambda^{-1} \big[   \frac 12 (\p_x u^h)^2     \big]
\\
&\le
\frac 12
  \|
(\p_x u^h)^2
\|_{L^2(\rr)}
\\
&\lesssim
   \|
\p_xu^h 
  \|_{L^\infty(\rr)}
  \|
\p_xu^h 
\|_{L^2(\rr)}
\\
&\lesssim
   \|
\p_xu^h 
  \|_{L^\infty(\rr)}
  \|
u^h 
\|_{H^1(\rr)}
\\
&\lesssim
   \lambda^{ -\frac 12 \delta - s +1}
\cdot
  \lambda^{-s +1},
 \end{split}
 \end{equation*}
 which gives
 \begin{equation}
 \label{H1-est-F6}
   \|F_6(t)\|_{H^1(\rr)}
\lesssim
 \lambda^{ -2s - \frac 12 \delta +2},
 \qquad 
 \lambda>>1.
 \end{equation}

 Collecting all error estimates together gives the following proposition.

%
 
%
%
%
%
%
%
%
%
%
%
\begin{proposition}
   \label{high-s-error-estimate-prop}
   Let $s>1$ and $1<\delta<2$. Then,
  for  $\omega$ bounded 
   and $ \lambda>>1$ we have that
   \begin{equation}
   \label{F-error-estimate-2}
   \|
   F(t)
   \|_{H^1(\rr)}
 \lesssim
 \lambda^{-r_s},
 \qquad
 \text{  for }
 \,
 \lambda >>1,
    \end{equation}
with
\begin{equation}
\label{decay-exponent}
 r_s
 \doteq
\big(s -  \frac 12\delta\big)>0,
\qquad
\text{ if }
\, 
s> \frac 12\delta.
\end{equation}
\end{proposition}
%
%
%

%
%
%
%
%

\section{  Estimating the difference between approximate and actual  solutions} 
\setcounter{equation}{0}
Let  $u_{\omega, \lambda}(x, t)$ be the solution to CH equation 
with initial data  the value of the approximate solution  $u^{\omega, \lambda}(x, t)$
at time zero.  That is,  $u_{\omega, \lambda}(x, t)$  solves  the Cauchy  problem
\begin{equation} 
\label{CH-with-appox-data} 
\partial_t u_{\omega, \lambda}
+
 u_{\omega, \lambda}
 \partial_xu_{\omega, \lambda}
 +
\Lambda^{-1}
 \Big[
 u_{\omega, \lambda}^2 + \frac{1}{2}(\partial_x  u_{\omega, \lambda})^2
\Big]
= 0, 
\ x \in  \mathbb{R},  \;\;  \ t \in \mathbb{R},
\end{equation}
\begin{equation} 
\label{CH-appox-data} 
u_{\omega, \lambda}(x, 0) 
=
u^{\omega, \lambda}(x, 0) 
=
\omega  \lambda^{-1}
\tilde{\varphi}(\frac{x}{\lambda^{\delta}})
+
  \lambda^{-\delta/2  -  s}
  \varphi (\frac{x}{\lambda^{\delta}})
  \cos (\lambda x).
\end{equation}
Note that  $u^{\omega, \lambda}(0)$ is in $H^s(\rr)$, $s\ge 0$, and 
   \begin{equation} 
\label{approx-sln-at-0-est} 
\|
u^{\omega, \lambda}(0)
\|_{H^s}
\le
\|
u_\ell(0)
\|_{H^{s}}
+
\|u^h(0)
\|_{H^s}
\lesssim 
\lambda^{-1+\frac 12 \delta} +1.
\end{equation}
Therefore, if $s>3/2$ then
  using  Theorem  \ref{CH-wp} and Proposition \ref{Lifespan-u-size}
 we see  that  for any  $\omega$ in a bounded set and  $\lambda>>1$
 the Cauchy problem \eqref{CH-with-appox-data}--\eqref{CH-appox-data} 
   has a unique solution $u_{\omega, \lambda}$ in  $C([0, T]; H^s(\rr))$ with
   \begin{equation} 
\label{T-indep-of-lambda} 
     T
        \gtrsim 
        \frac{1}{ 
     \|
u^{\omega, \lambda}(0)
  \|_{H^s(\rr)}}
    \gtrsim 
        \frac{1}{ 
   \lambda^{-1+\delta/2} +1} 
    \gtrsim   1.
 \end{equation}
   In fact,    $u_{\omega, \lambda}(t)$ is in $C^\infty$  for each $t\in [0, T]$. 
   
 %
%

To estimate the  difference between approximate and actual solutions 
we form the differential equation which it satisfies. So, if we let
     \begin{equation}
   \label{CH-difference}
v
=
u^{\omega, \lambda} -  u_{\omega, \lambda},
 \end{equation}
then a straightforward computation shows that $v$ satisfies the
Cauchy problem 
 \begin{equation} 
\label{CH-difference-eqn} 
\p_t v
-
 v\p_xv
 +
u^{\omega, \lambda} \p_xv
 +
 \p_xu^{\omega, \lambda} v
-
 \Lambda^{-1}
 \Big[
 v^2 
  +
  \frac{1}{2}(\partial_x  v)^2
 -
 2u^{\omega, \lambda} v
  -
  \p_xu^{\omega, \lambda} \p_xv
\Big]
= F(x, t),
\end{equation}
\begin{equation} 
\label{CH-difference-data} v(x, 0) = 0,
\ x \in  \mathbb{R},  \;\;  \ t \in \mathbb{R}.
\end{equation}
where  $F$ is defined by
\begin{equation} 
\label{Burgers-F} 
F
\doteq
 \p_tu^{\omega, \lambda}
+
u^{\omega, \lambda}
\p_xu^{\omega, \lambda}
+
\Lambda^{-1}
   [(u^{\omega, \lambda})^2
+
\frac 12 (\p_xu^{\omega, \lambda})^2],
\end{equation}
and which it has been shown to satisfy the $H^1$-estimate \eqref{F-error-estimate-2}.

 %
 %
%
%
\begin{lemma}
\label{CH-differ-H1-est-lem} 
Let $1<\delta<2$. 
If    $s>3/2$  then
  \begin{equation} 
  \|
v(t)
\|_{H^1(\rr)}
\doteq
\label{differ-H1-est} 
\|
u^{\omega, \lambda}(t) 
- 
 u_{\omega, \lambda}(t)
\|_{H^1(\rr)}
\lesssim 
 \lambda^{-r_s}, 
 \quad
  0 \le t \le T,
\end{equation}
where $r_s=s- \delta/2>0$  (see \eqref{decay-exponent}).
\end{lemma}
 
 %
%

%
{\bf Proof.}   We have
\begin{equation} 
\label{t-deriv-of-H1-v}
\frac 12
\frac{d}{dt}
\|
v(t)
\|_{H^1(\rr)}^2
=
\int_{\rr} 
\big[
v\partial_tv
+
\partial_xv \partial_x\partial_tv
\big]
dx
\end{equation}
Applying to both sides of  \eqref{CH-difference-eqn} the operator $(1-\partial_x^2)$
and solving for $\partial_t v$ we obtain
\begin{equation} 
\label{t-deriv-of-v}
\begin{split} 
\partial_t v
&=
(1-\partial_x^2)  F
\\
&
-
  (1-\partial_x^2)   \big[
u^{\omega, \lambda} \p_xv
 +
 \p_xu^{\omega, \lambda} v
   \big]
 \\
&
 -
 \partial_x
 \big[
2u^{\omega, \lambda} v
  +
  \p_xu^{\omega, \lambda} \p_xv
  \big]
  \\
  &
  +
   3v\partial_x v
-
2
 \partial_xv\partial_x^2 v
 -
 v\partial_x^3 v
 +
 \partial_t \partial_x^2v
\end{split}
\end{equation}  
Substituting $\partial_tv$ from \eqref{t-deriv-of-v}
to \eqref{t-deriv-of-H1-v}  we get
\begin{equation} 
\label{deriv-H1-norm-v-sq-1}
\begin{split} 
\frac 12
\frac{d}{dt}
\|
v(t)
\|_{H^1(\rr)}^2
&
=
\int_{\rr} v (1-\partial_x^2)  F dx
\\
&
-
\int_{\rr}
v  (1-\partial_x^2)   \big[
u^{\omega, \lambda} \p_xv
 +
 \p_xu^{\omega, \lambda} v
 \big]
 dx
 \\
&
 -
\int_{\rr}
v
 \partial_x
 \big[
2u^{\omega, \lambda} v
  +
  \p_xu^{\omega, \lambda} \p_xv
\big]
dx
\\
&
 +
\int_{\rr}
 \big[
v
   (3v\partial_x v
-
 2\partial_xv\partial_x^2 v
 -
 v\partial_x^3 v
 +
 \partial_t \partial_x^2v)
 +
 \partial_xv \partial_x\partial_tv
\big]
dx.
\end{split}
\end{equation}
Noting that the last integral can be rewritten as 
$$
\int_{\rr} 
\Big[ 
\partial_x\big(v^3\big)
-
\partial_x\big(v^2 \partial_x^2v\big)
+
\partial_x\big(v\partial_t\partial_xv\big)
 \Big] \, dx
  =
0,
$$
which is a property special to CH, 
we see that 
equation \eqref{deriv-H1-norm-v-sq-1}  takes the form 
\begin{equation} 
\label{deriv-H1-norm-v-sq}
\begin{split}
\frac 12 
\frac{d}{dt}
\|
v(t)
\|_{H^1(\rr)}^2
&
=
\int_{\rr} 
v (1-\partial_x^2)  F dx
\\
&
-
\int_{\rr}
v  (1-\partial_x^2)   \big[
u^{\omega, \lambda} \p_xv
 +
 \p_xu^{\omega, \lambda} v
 \big]
 dx
 \\
&
 -
\int_{\rr}
v
 \partial_x
 \big[
2u^{\omega, \lambda} v
  +
  \p_xu^{\omega, \lambda} \p_xv
\big]
dx.
\end{split}
\end{equation}
 Integrating by parts and  applying the Cauchy-Schwarz inequality, 
 we estimate the three integrals  in  the right-hand side
 of   \eqref{deriv-H1-norm-v-sq} as follows.
 For the first  integral we have
 \begin{equation} 
\label{RHS-1}
\Big|
 \int_{\rr} 
v (1-\partial_x^2)  F dx
\Big|
=
\Big|
 \int_{\rr} 
 [ 
  v F
 +
\p_xv\p_xF 
]
dx
\Big|
\le
\|
F(t)
\|_{H^1(\rr)}
\|
v(t)
\|_{H^1(\rr)}.
\end{equation}
 Also, for the third integral we have
 \begin{equation} 
\label{RHS-3}
\begin{split}
\Big|
\int_{\rr}
v
 \partial_x
 \big[
2u^{\omega, \lambda} v
  &+
  \p_xu^{\omega, \lambda} \p_xv
\big]
dx
\Big|
=
\Big|
\int_{\rr}
\p_xv
 \big[
2u^{\omega, \lambda} v
  +
  \p_xu^{\omega, \lambda} \p_xv
\big]
dx
\Big|
\\
&
\le
2
\Big(
\|
u^{\omega, \lambda}(t)
\|_{L^\infty(\rr)}
+
\|
\p_xu^{\omega, \lambda}(t)
\|_{L^\infty(\rr)}
\Big)
\|
v(t)
\|_{H^1(\rr)}^2.
\end{split}
\end{equation}
Integrating by parts, we write the second  integral  in the form
 \begin{equation} 
\label{RHS-3}
\begin{split}
\int_{\rr}
v  (1-\partial_x^2)   \big[
u^{\omega, \lambda} \p_xv
 &+
 \p_xu^{\omega, \lambda} v
 \big]
 dx
=
\int_{\rr}
v    \big[
u^{\omega, \lambda} \p_xv 
 +
 \p_xu^{\omega, \lambda} v
 \big]  dx
\\
&+
 \int_{\rr}
\p_x v   \p_x \big[
u^{\omega, \lambda} \p_xv \big] dx
+
 \int_{\rr}
 \p_x v   \p_x 
\big[
 \p_xu^{\omega, \lambda} v
 \big]
  dx
\end{split}
\end{equation}
and estimate its first part by 
 \begin{equation} 
\label{RHS-2.1}
\Big|
\int_{\rr}
v    \big[
u^{\omega, \lambda} \p_xv 
 +
 \p_xu^{\omega, \lambda} v
 \big]
 dx
 \Big|
\le
\Big(
\|
u^{\omega, \lambda}(t)
\|_{L^\infty(\rr)}
+
\|
\p_xu^{\omega, \lambda}(t)
\|_{L^\infty(\rr)}
\Big)
\|
v(t)
\|_{H^1(\rr)}^2.
\end{equation}
Its second part we can be written as 
 \begin{equation*} 
 \label{RHS-2.2}
 \int_{\rr}
\p_x v   \p_x \big[
u^{\omega, \lambda} \p_xv \big] dx
=
\int_{\rr}
 \big[
 \frac 12 u^{\omega, \lambda} \p_x(\p_xv)^2
 +
\p_x u^{\omega, \lambda} (\p_xv)^2
 \big] dx
 =
 \frac 12
 \int_{\rr}
 \p_xu^{\omega, \lambda} (\p_xv)
 \big] dx,
\end{equation*}
which gives that
 \begin{equation} 
 \label{RHS-2.2}
 \Big|
 \int_{\rr}
\p_x v   \p_x \big[
u^{\omega, \lambda} \p_xv \big] dx
 \Big|
\le
\|
\p_xu^{\omega, \lambda}(t)
\|_{L^\infty(\rr)}
\|
v(t)
\|_{H^1(\rr)}^2.
\end{equation}
Finally,  writing the  last part as follows
 \begin{equation*} 
 \label{RHS-2.3}
 \int_{\rr}
 \p_x v   \p_x 
\big[
 \p_xu^{\omega, \lambda} v
 \big]
  dx
=
\int_{\rr}
 \big[
\p_xu^{\omega, \lambda} (\p_xv)^2
 +
\p_x^2 u^{\omega, \lambda}  v  \p_xv
 \big] dx
\end{equation*}
we see that it can be estimated   as follow
 \begin{equation} 
 \label{RHS-2.3}
 \Big|
  \int_{\rr}
 \p_x v   \p_x 
\big[
 \p_xu^{\omega, \lambda} v
 \big]
  dx
 \Big|
\le
\Big(
\|
\p_xu^{\omega, \lambda}(t)
\|_{L^\infty(\rr)}
+
\|
\p_x^2u^{\omega, \lambda}(t)
\|_{L^\infty(\rr)}
\Big)
\|
v(t)
\|_{H^1(\rr)}^2.
\end{equation}
Combining  the above estimates  gives
  \begin{equation} 
\label{deriv-H1-norm-v-sq-2}
\begin{split} 
\frac 12
\frac{d}{dt}
&\|
v(t)
\|_{H^1(\rr)}^2
\lesssim 
\|
F(t)
\|_{H^1(\rr)}
\|
v(t)
\|_{H^1(\rr)}
+
\\
&
+
\Big(
\|
u^{\omega, \lambda}(t)
\|_{L^\infty(\rr)}
+
\|
\p_xu^{\omega, \lambda}(t)
\|_{L^\infty(\rr)}
+
\|
\p_x^2u^{\omega, \lambda}(t)
\|_{L^\infty(\rr)}
\Big)
\|
v(t)
\|_{H^1(\rr)}^2.
\end{split}
\end{equation}
  From   \eqref{high-frequency-approx-sln}  we have
\begin{equation} 
\label{CH-xx-deriv}
\begin{split}
\p_x^2u^h
&= 
  \lambda^{-\frac 52 \delta  -  s}
  \p_x^2\varphi (\frac{x}{\lambda^{\delta}})
  \cos (\lambda x-\omega t)
  \\
  &
  - 2
  \lambda^{-\frac 32 \delta  -  s+1}
  \p_x\varphi (\frac{x}{\lambda^{\delta}})
  \sin (\lambda x-\omega t)
   - 2
  \lambda^{-\frac 12 \delta  -  s+2}
  \varphi (\frac{x}{\lambda^{\delta}})
  \cos (\lambda x-\omega t).
  \end{split}
 \end{equation}
 so  that   
  \begin{equation} 
\label{u-h-ap-x-sup}  
\|
u^h(t)
\|_{L^\infty(\rr)}
+
\|
\p_xu^h(t)
\|_{L^\infty(\rr)}
+
\|
\p_x^2u^h(t)
\|_{L^\infty(\rr)}
\lesssim
\lambda^{-({\frac 12 \delta  + s-2})}.
\end{equation}
For $u_\ell$ we have 
 \begin{equation} 
\label{u-ell-x-sup}  
\|
u_\ell(t)
\|_{L^\infty(\rr)}
+
\|
\p_xu_\ell(t)
\|_{L^\infty(\rr)}
+
\|
\p_x^2u_\ell(t)
\|_{L^\infty(\rr)}
\lesssim
\|
u_\ell(t)
\|_{H^3(\rr)}
\lesssim
\lambda^{-(1-\frac 12 \delta)}.
\end{equation}
Therefore
  \begin{equation} 
\label{u-ap-x-sup}  
\|
u^{\omega, \lambda}(t)
\|_{L^\infty(\rr)}
+
\|
\p_xu^{\omega, \lambda}(t)
\|_{L^\infty(\rr)}
+
\|
\p_x^2u^{\omega, \lambda}(t)
\|_{L^\infty(\rr)}
\lesssim
\lambda^{-\rho_s},
\end{equation}
where
  \begin{equation} 
\label{rho-def}  
\rho_s
\doteq
\min\{1-\frac 12 \delta,  \frac 12 \delta  + s-2
\}>0,
\end{equation}
for any  any \textcolor{Red}{ $s>3/2$} if   $\delta$ is chosen appropriately
in the interval $(1, 2)$.

Using \eqref{u-ap-x-sup}   and the $H^1$- estimate    \eqref{F-error-estimate-2}
 for the  error $F$, from \eqref{deriv-H1-norm-v-sq-2}  we get
\begin{equation*} 
\label{deriv-H1-norm-v-sq-3}
\frac 12
\frac{d}{dt}
\|
v(t)
\|_{H^1(\rr)}^2
\lesssim 
 \lambda^{-\rho_s}
\|
v(t)
\|_{H^1(\rr)}^2
+
 \lambda^{-r_s}
\|
v(t)
\|_{H^1(\rr)},
\end{equation*}
which gives  the differential inequality 
\begin{equation} 
\label{deriv-H1-norm-v-sq-3}
\frac{d}{dt}
\|
v(t)
\|_{H^1(\rr)}
\lesssim 
 \lambda^{-\rho_s}
\|
v(t)
\|_{H^1(\rr)}
+
 \lambda^{-r_s}.
\end{equation}
Since $ \|v(0)\|_{H^1(\rr)}=0$ and for 
 $s>1$ we can choose  $\delta$   such  that  $\rho_s \ge 0$
from \eqref{deriv-H1-norm-v-sq-3} and Gronwall's inequality  we obtain that
  \begin{equation} 
\label{v-differ-est} 
\|
v(t)
\|_{H^1(\rr)}
\lesssim 
 \lambda^{-r_s}, 
 \quad
  0 \le t \le T,
\end{equation}
which concludes  the proof of the lemma. \, $\square$
%

%
%
%
%
%
%

\section{  Non-uniform dependence in $H^s(\rr)$ for $s>3/2$} 
\setcounter{equation}{0}
 Next we  shall  prove non-uniform dependence for CH
 by taking  advantage of the  information provided 
 by  Theorem \ref{CH-wp}   and Proposition \ref{Lifespan-u-size},
and  the $H^1$- estimate \eqref{differ-H1-est} on
 the difference between  approximate solutions and solutions
 with   same initial data. 
 
 For this, let  $u_{1, \lambda}(x, t)$ and   $u_{-1, \lambda}(x, t)$ be the unique solutions 
   to the  the Cauchy problem  
   \eqref{CH-with-appox-data}--\eqref{CH-appox-data} 
   with initial data  
   $u^{1, \lambda}(x, 0)$ and   $u^{-1, \lambda}(x, 0)$ correspondingly. 
   By Theorem \ref{CH-wp}  these solutions belong in  $C([0, T]; H^s(\rr))$.
   Recall, using  Proposition \ref{Lifespan-u-size} we proved
 estimate \eqref{T-indep-of-lambda}  which  says that 
 $T$ is independent of  $\lambda>>1$. Also,  for $s>3/2$,
  using estimate \eqref{u-u0-Hs-bound},  we have
     \begin{equation}
   \label{CH-slns-Hs-norm}
\|
u_{\pm 1, \lambda}(t)
  \|_{H^s(\rr)}
  \lesssim
  \|
u^{\pm 1, \lambda}(0)
  \|_{H^s(\rr)},
  \quad 0\le t\le T.
   \end{equation}
   %
   %
   %
   %
   %
   %
Furthermore,  since our $s$-dependent  initial data 
$u^{\pm, \lambda}(0)$  belong to every Sobolev space
they do belong to $H^{[s]+2}(\rr)$.
Since $s>3/2$ by the argument in the last  remark 
of section 2 we obtain a companion estimate 
 to  \eqref{CH-slns-Hs-norm} 
     \begin{equation}
   \label{CH-slns-Hk-norm}
\|
u_{\pm 1, \lambda}(t)
  \|_{H^{[s]+2}(\rr)}
  \lesssim
  \|
u^{\pm 1, \lambda}(0)
  \|_{H^{[s]+2}(\rr)}, \quad 0\le t\le T.
   \end{equation}
Now let $k=[s]+2$.
  If   $\lambda$ is large enough   then  from
 \eqref{Hs-norm-of-ap-sl-2} and  \eqref{Hs-norm-of u-ell-t-est}
 we have
     \begin{equation*}
   \label{ap-data-size}
   \begin{split}
 \|
u^{\pm 1, \lambda}(t)
  \|_{H^k(\rr)}
 & \le  
 \|
u_{\ell, \pm 1, \lambda}(t)
  \|_{H^k(\rr)}
  +
   \lambda^{-\frac{1}{2}\delta-s}
  \|
  \varphi(\frac{x}{\lambda^\delta})  \cos (\lambda x-\lambda t)
  \|_{H^k(\rr)}
  \\
  &
\lesssim
  \lambda^{-1+\frac 12 \delta}
  +
     \lambda^{k-s} \cdot \lambda^{-\frac{1}{2}\delta-k}
  \|
  \varphi(\frac{x}{\lambda^\delta})  \cos (\lambda x-\lambda t)
  \|_{H^k(\rr)}
\\
  &
\lesssim
  \lambda^{-1+\frac 12 \delta}
  +
     \lambda^{k-s}
 \|
  \varphi
  \|_{L^2(\rr)},
    \end{split}
 \end{equation*}
which gives
     \begin{equation}
   \label{ap-data-Hk-size}
 \|
u^{\pm 1, \lambda}(t)
  \|_{H^k(\rr)}
  \lesssim
   \lambda^{k-s},
   \quad
   \text{hence by   \eqref{CH-slns-Hk-norm} }
      \quad
   \|
u_{\pm 1, \lambda}(t)
  \|_{H^k(\rr)}
    \lesssim
   \lambda^{k-s}.
   \quad 
    \end{equation}
 Therefore, from  \eqref{ap-data-Hk-size}  we 
 obtain the following  estimate for the $H^k$-norm of the difference  of
 $u_{\pm 1, \lambda}$ and  $u_{\pm 1, \lambda}$
     \begin{equation}
   \label{CH-differ-Hk-norm}
\|
u^{\pm 1, \lambda}(t)
-
u_{\pm 1, \lambda}(t)
  \|_{H^k(\rr)}
    \lesssim
     \lambda^{k-s},
   \quad
  0 \le t \le T.
   \end{equation}
   Applying  \eqref{differ-H1-est}  with our particular choice of $\omega=\pm 1$ we have
  \begin{equation} 
\label{CH-differ-H1-est} 
\|
u^{\pm 1, \lambda}(t)
-
u_{\pm 1, \lambda}(t)
  \|_{H^1(\rr)}
\lesssim 
 \lambda^{-r_s }, 
 \quad
  0 \le t \le T.
\end{equation}
Now, applying  the  interpolation inequality
  \begin{equation*} 
\label{Hs-interpolation} 
\|
\psi
  \|_{H^s(\rr)}
\le
\|
\psi
  \|_{H^{s_1}(\rr)}^{(s_2-s)/(s_2-s_1)}
  \|
\psi
  \|_{H^{s_2}(\rr)}^{(s-s_1)/(s_2-s_1)}
\end{equation*}
with $s_1=1$ and $s_2=[s] +2=k$ and 
 using estimates \eqref{CH-differ-H1-est} and  \eqref{CH-differ-Hk-norm}
gives
  \begin{equation} 
\label{CH-differ-Hs-est-interpo} 
\begin{split}
\|
u^{\pm 1, \lambda}(t)
-
u_{\pm 1, \lambda}(t)
  \|_{H^s(\rr)}
&\le
\|
u^{\pm 1, \lambda}(t)
-
u_{\pm 1, \lambda}(t)
  \|_{H^1(\rr)}^{(k-s)/(k-1)}
  \\
  &\cdot
   \|
  u^{\pm 1, \lambda}(t)
-
u_{\pm 1, \lambda}(t)
  \|_{H^k(\rr)}^{(s-1)/(k-1)}
  \\
  &
  \lesssim
 \lambda^{(-r_s)[(k-s)/(k-1)]} 
  \lambda^{(k-s)[(s-1)/(k-1)]}
   \\
  &
  \lesssim
  \lambda^{-(r_s-s+1)[(k-s)/(k-1)]}.
    \end{split}
\end{equation}
From the last inequality we obtain that
  \begin{equation} 
\label{CH-differ-Hs-dacay-est} 
\|
u^{\pm 1, \lambda}(t)
-
u_{\pm 1, \lambda}(t)
  \|_{H^s(\rr)}
   \lesssim
 \lambda^{-\varepsilon_s}, 
 \,\,\,
  0 \le t \le T,
  \end{equation}
  where  $\varepsilon_s$ is given by
\begin{equation}
\label{epsilon-decay-exponent}
 \varepsilon_s
 = (1 - \frac 12 \delta)/(s+2).
 \end{equation}
Note  that
\begin{equation}
\label{epsilon-s-positive}
 \varepsilon_s>0,
 \quad
 \text{for}
 \quad
s>1.
\end{equation}

  Next, we shall use  estimate \eqref{CH-differ-Hs-dacay-est} to prove non-uniform dependence  when 
$s>3/2$.
%
%

  %
%

%
\medno
{\bf Behavior at time zero.}  Since  $\delta<2$, at $t=0$ we have
  \begin{equation} 
\label{CH-slns-differ-t-0} 
\begin{split}
\|
u_{1, \lambda}(0)
-
u_{-1, \lambda}(0)
  \|_{H^s(\rr)}
&=
\|
2  \lambda^{-1}
\tilde{\varphi}(\frac{x}{\lambda^{\delta}})
  \|_{H^s(\rr)}
  \\
&\le
2
  \lambda^{-1+\frac 12 \delta}
 \|
\tilde{\varphi}
  \|_{H^s(\rr)}
  \longrightarrow
  0
  \,\,
  \text{as}
  \,\,
  \lambda \to \infty.
  \end{split}
  \end{equation}
  %
  %
%

%
\nin
{\bf Behavior at time  $t>0$.}  Then, we write
  \begin{equation} 
  \label{CH-slns-differ-t-pos}
\begin{split}
\|
u_{1, \lambda}(t)
-
u_{- 1, \lambda}(t)
  \|_{H^s(\rr)}
  &
\ge
\|
u^{1, \lambda}(t)
-
u^{- 1, \lambda}(t)
  \|_{H^s(\rr)}
  \\
  &
  -
  \|
u^{1, \lambda}(t)
-
u_{1, \lambda}(t)
  \|_{H^s(\rr)}
  \\
  &
  -
  \|
u^{-1, \lambda}(t)
-
u_{-1, \lambda}(t)
  \|_{H^s(\rr)}
   \end{split}
\end{equation}
Using estimate \eqref{CH-differ-Hs-dacay-est} for the last two terms 
in \eqref{CH-slns-differ-t-pos} we obtain
  \begin{equation} 
  \label{CH-slns-differ-t-pos-est}
\|
u_{1, \lambda}(t)
-
u_{- 1, \lambda}(t)
  \|_{H^s(\rr)}
\ge
\|
u^{1, \lambda}(t)
-
u^{- 1, \lambda}(t)
  \|_{H^s(\rr)}
-
c
 \lambda^{-\varepsilon_s}.
\end{equation}
In  \eqref{CH-slns-differ-t-pos-est} letting $\lambda$ go to $\infty$ gives
  \begin{equation} 
  \label{CH-slns-to-ap-est}
  \liminf_{\lambda\to\infty}
\|
u_{1, \lambda}(t)
-
u_{- 1, \lambda}(t)
  \|_{H^s(\rr)}
\ge
  \liminf_{\lambda\to\infty}
\|
u^{1, \lambda}(t)
-
u^{- 1, \lambda}(t)
  \|_{H^s(\rr)}.
\end{equation}
Inequality  \eqref{CH-slns-to-ap-est} is a key estimate since it reduces  finding  a lower positive bound 
for the difference of the  {\bf uknown} solution sequences to
 finding  a lower positive bound 
for the difference of the known approximate solution sequences.
Using the identity
$$
\cos \alpha -\cos \beta
=
-2
\sin(\frac{\alpha + \beta}{2})
\sin(\frac{\alpha - \beta}{2})
$$
gives
$$
u^{1, \lambda}(t)
-
u^{- 1, \lambda}(t)
=
u_{\ell, 1, \lambda}(t)
-
u_{\ell,  -1, \lambda}(t)
+
2
  \lambda^{-\frac 12 \delta  -  s}
  \varphi (\frac{x}{\lambda^{\delta}})
  \sin (\lambda x) \sin t.
  $$
  Therefore
  \begin{equation} 
  \label{B--ap-below-est-1}
  \begin{split}
\|
u^{1, \lambda}(t)
-
u^{- 1, \lambda}(t)
  \|_{H^s(\rr)}
  &
  \ge
  2
  \lambda^{-\frac 12 \delta  -  s}
  \|
  \varphi (\frac{x}{\lambda^{\delta}})
  \sin (\lambda x)
  \|_{H^s(\rr)}
   |\sin t|
   \\
   & -
   \|
   u_{\ell, 1, \lambda}(t)
  \|_{H^s(\rr)}
  -
   \|
   u_{\ell, -1, \lambda}(t)
  \|_{H^s(\rr)}
    \\
  &
 \gtrsim
  2
  \lambda^{-\frac 12 \delta  -  s}
  \|
  \varphi (\frac{x}{\lambda^{\delta}})
  \sin (\lambda x)
  \|_{H^s(\rr)}
   |\sin t|
   -
  \lambda^{-1+\frac 12 \delta}.
    \end{split}
\end{equation}
Now    letting $\lambda$ go to $\infty$,  \eqref{B--ap-below-est-1}  gives
  \begin{equation} 
  \label{CH-ap-below-est}
  \liminf_{\lambda\to\infty}
\|
u^{1, \lambda}(t)
-
u^{- 1, \lambda}(t)
  \|_{H^s(\rr)}
 \gtrsim
  \|
\varphi
  \|_{L^2(\rr)}
    |\sin t|.
\end{equation}
Combining  \eqref{CH-slns-to-ap-est} and  \eqref{CH-ap-below-est} gives
  \begin{equation} 
  \label{CH-slns-below-est-fin}
  \liminf_{\lambda\to\infty}
\|
u_{1, \lambda}(t)
-
u_{- 1, \lambda}(t)
  \|_{H^s(\rr)}
 \gtrsim
  \|
\varphi
  \|_{L^2(\rr)}
    |\sin t  |,
\end{equation}
which proves  Theorem \ref{CH-non-unif-dependence}.

\vskip0.1in
\noindent
{\em{ \bf ACKNOWLEDGEMENTS}}. The first author thanks the Department
of Mathematics of the University of Chicago for the hospitality
during his stay there in the Fall  of  2007.
The second author acknowledges  partial support from the  NSF.


\vskip0.1in

\begin{minipage}[b]{6 cm} {\bf A. Alexandrou Himonas}\\ Department of
Mathematics \\
     University of Notre Dame\\ Notre Dame, IN 46556\\
     E-mail: {\it himonas.1$@$nd.edu}
\end{minipage}
\hfill
\begin{minipage}[b]{7 cm}
{\bf Carlos Kenig}\\
Department of Mathematics\\
The University of Chicago\\
5734 S. University Avenue\\
Chicago, Illinois 60637\\
 E-mail: {\it cek$@$math.uchicago.edu}
\end{minipage}

\end{document}